\documentclass[a4paper,11pt]{article}

\usepackage[utf8]{inputenc}
\usepackage[english,activeacute]{babel}
\usepackage{amsmath,amsthm,amsfonts,amssymb}
\usepackage{mathpartir}
\usepackage{ stmaryrd }
\usepackage{graphics}
\usepackage{enumerate}
\usepackage{mathrsfs}
\usepackage{rotating}
\usepackage{hyperref}
\usepackage{xcolor}
\hypersetup{
    colorlinks,
    linkcolor={red!50!black},
    citecolor={blue!50!black},
    urlcolor={blue!80!black}
}
\input{xy}
\usepackage[all]{xy}
\usepackage{tikz}
\usepackage{tikz-cd}

\newtheorem{thm}{Theorem}[section]
\newtheorem{cor}[thm]{Corollary}
\newtheorem{lemma}[thm]{Lemma}

\newtheorem{rmk}[thm]{Remark}
\theoremstyle{plain}
\newtheorem{defs}[thm]{Definition}
\theoremstyle{remark}
\newtheorem{exm}[thm]{Example}

\input{diagxy}
\xyoption{curve}

\newbox\anglebox % large pullback angle
\setbox\anglebox=\hbox{\xy \POS(75,0)\ar@{-} (0,0) \ar@{-} (75,75)\endxy}
 
\newbox\angleboxr % reverse large pullback angle
\setbox\angleboxr=\hbox{\xy \POS(0,0)\ar@{-} (0,75) \ar@{-} (75,0)\endxy}
 
\newbox\sanglebox % small pullback angle
\setbox\sanglebox=\hbox{\xy \POS(50,0)\ar@{-} (0,0) \ar@{-} (50,50)\endxy}
 
\newbox\sangleboxr % small reverse pullback angle
\setbox\sangleboxr=\hbox{\xy \POS(0,0)\ar@{-} (0,50) \ar@{-} (50,0)\endxy}
 
\newbox\sangleboxf % small flipped pullback angle
\setbox\sangleboxf=\hbox{\xy \POS(50,50)\ar@{-} (50,0) \ar@{-} (0,50)\endxy}
 
\newbox\angleboxf % flipped pullback angle
\setbox\angleboxf=\hbox{\xy \POS(75,75)\ar@{-} (75,0) \ar@{-} (0,75)\endxy}
 
\newbox\sangleboxfr % small flipped reverse pullback angle
\setbox\sangleboxfr=\hbox{\xy \POS(0,50)\ar@{-} (50,50) \ar@{-} (0,0)\endxy}
 
\newbox\angleboxfr % small flipped reverse pullback angle
\setbox\angleboxfr=\hbox{\xy \POS(0,75)\ar@{-} (75,75) \ar@{-} (0,0)\endxy}

\newcommand{\Sets}{\ensuremath{\mathbf{Set}}}

\title{A proof of Shelah's eventual categoricity conjecture and an extension to accessible categories with directed colimits}
\author{Christian Esp\'indola}

\begin{document}
\date{}
\maketitle

\begin{abstract}
We provide a proof, in $ZFC$, of Shelah's eventual categoricity conjecture for abstract elementary classes (AEC's). Moreover, assuming in addition the Singular Cardinal Hypothesis ($SCH$), we prove a direct generalization to the more general context of accessible categories with directed colimits. If $\mathcal{K}$ is such a category, we show that there is a cardinal $\mu$ such that if $\mathcal{K}$ is $\lambda$-categorical for some $\lambda \geq \mu$ (i.e., it has only one object of internal size $\lambda$ up to isomorphism), then $\mathcal{K}$ is eventually categorical (i.e., it is $\lambda'$-categorical for every $\lambda' \geq \mu$). When considering cardinalities of models of infinitary theories $\mathbb{T}$ of $\mathcal{L}_{\kappa, \theta}$ that axiomatize $\mathcal{K}$, the result implies, under $SCH$, the following infinitary version of Morley's categoricity theorem: let $S$ be the class of cardinals $\lambda$ which are of cofinality at least $\theta$ but are not successors of cardinals of cofinality less than $\theta$. Then, if $\mathbb{T}$ is a $\mathcal{L}_{\kappa, \theta}$ theory whose models have directed colimits and it is $\lambda$-categorical for some $\lambda \geq \mu$ in $S$, then it is $\lambda'$-categorical for every $\lambda' \geq \mu$ in $S$; moreover, we also exhibit an example that shows that the exceptions in the class $S$ are needed. Along the way we also prove Grossberg conjecture, according to which categoricity in a high enough cardinal implies eventual amalgamation. We establish this result in AEC's and, assuming in addition $SCH$, in the more general context of accessible categories whose morphisms are monomorphisms.\footnote{This research has been partially supported through the grant 19-00902S from the Grant Agency of the Czech Republic.}
\end{abstract}

\section{Introduction}

  This paper belongs to the classification theory program initiated by Shelah, intended to determine dividing lines in the behaviour of the categoricity spectrum of theories. Continuing the work of Morley, who established that categoricity of a countable finitary theory in one uncountable cardinal entailed categoricity in all uncountable cardinals, Shelah intended to generalize this setup, to uncountable finitary theories and even to finite quantifier infinitary theories. He later pushed the classification to more general infinitary theories introducing abstract elementary classes. The main open problem in the area, Shelah's eventual categoricity conjecture, asserts that for any abstract elementary class (AEC) there is a cardinal $\kappa$ such that if the AEC is categorical in some $\lambda \geq \kappa$, it is categorical in \emph{all} $\lambda' \geq \kappa$. This general conjecture was stated in \cite{shelah}, while the version for the particular case of sentences in $\mathcal{L}_{\omega_1, \omega}$ was conjectured around 1977. Several approximations are known (see, e.g., the account presented in the introduction of \cite{sv3}), though the full conjecture has been open until now. When there is a proper class of strongly compact cardinals or when amalgamation and a suitable weakening of the Generalized Continuum Hypothesis ($GCH$) hold, the conjecture was proven to be true in \cite{shva}. We will provide here a full proof of the conjecture, in $ZFC$, and show as well that eventual amalgamation in fact follows from categoricity in a high enough cardinal, thus proving a conjecture of Grossberg from 1986. Moreover, our topos-theoretic argument is general enough so that, assuming in addition $SCH$, it applies to accessible categories with directed colimits, establishing these conjectures in a much wider variety of categories of models. This also provides a version of Morley's categoricity theorem for the corresponding infinitary theories $\mathbb{T}$ in $\mathcal{L}_{\kappa, \theta}$ whose models have directed colimits; if $S$ be the class of cardinals $\lambda$ which are of cofinality at least $\theta$ but are not successors of cardinals of cofinality less than $\theta$ then, under $SCH$, when $\mathbb{T}$ is $\lambda$-categorical for some $\lambda \geq \mu$ in $S$, it is $\lambda'$-categorical for every $\lambda' \geq \mu$ in $S$. The case $\theta=\omega$ is (a generalization of) Shelah's original conjecture, since in that case the use of $SCH$ can be eliminated\footnote{In fact, the precise amount of $SCH$ that is required to hold for infinitary theories $\mathbb{T}$ in $\mathcal{L}_{\kappa, \theta}$ is $SCH_{\theta, \geq 2^{<\theta}}$, defined as ``for all $\mu \geq 2^{<\theta}$ there is a set of cardinals $\lambda_i \leq \mu$ unbounded below $\mu$ such that, for each $i$, $\nu^{<\theta} \leq \lambda_i$ for all $\nu<\lambda_i$". This is discussed in Remark 2.3 of \cite{rlv}.}; also, the threshold cardinal we find is not explicit, but it is provably the least possible for which the conjecture is true. We prove as well that this infinitary version of Morley's theorem is sharp in the sense that we build examples of theories (see Example \ref{exm}) where the cardinals $\lambda$ which are not in the class $S$ fail in fact to have categoricity in the cardinality of the model, making our result best possible. 

  The main tool that we will use for these results is a completeness theorem for infinitary intuitionistic logics developed by the author in \cite{espindola} and \cite{espindolad}. These are theorems valid for theories in infinite quantifier languages, for which a new sound rule of inference has been devised. Even when restricting to the particular case of classical theories, these completeness theorems are sharper than known completeness theorems for infinitary classical logic (e.g., those developed by Karp in \cite{karp}), as unlike the latter, the former are able to be cast in the setting of categorical logic, making use of an infinitary generalization of syntactic categories that begun to be explored in the '90's for the infinitary regular fragment. As we will see, the categorical language is essential in this paper since the crucial arguments are topos-theoretic in nature. What lies below is the fact that any abstract elementary class is equivalent to the category of $\kappa$-points of some $\kappa$-topos (introduced in \cite{espindolad}, where they are called $\kappa$-geometric toposes). It is this precise observation which will allow for our generalization to accessible categories with directed colimits.

  Among the consequences of the completeness theorem that we will use, the omitting types theorem for infinite quantifier languages stands out. Whenever the category of models of an infinitary theory over an infinite quantifier language has directed colimits, an analogue of the well-known omitting types theorem for finitary logic (or $\mathcal{L}_{\omega_1, \omega}$) can be derived as a consequence of our completeness theorem. It was realized in the '70's that the omitting types theorem was an essential tool in model theory, having consequences for atomic and prime models, among other applications; this will be recast in our case to study infinitary theories which are categorical in some cardinal. The last ingredient of our proof is a precise characterisation of the $\kappa$-classifying topos for $\kappa$-saturated models, which will be very useful for deriving a downward categoricity transfer.

  Our approach avoids large cardinals or amalgamation since these are essentially, as we will see, compactness properties of the theories in question, and we will be able to derive them directly from a categoricity assumption by using instead an approach that relies on infinitary regular logic, introduced by Makkai in \cite{makkai}. Although categoricity is not enough to make our theory axiomatizable in infinitary regular logic, it can be axiomatized in a slightly bigger fragment, which can then be studied through an equivalent infinitary regular theory on a different signature, an instance of the process known as Morleyization. $\kappa$-regular logic is known to be $\kappa$-compact, as proven already by Makkai in \cite{makkai}. The reason is essentially that while compactness arguments usually involve \L{}o\'s theorem and the consequent use of $\kappa$-complete ultrafilters extending a given $\kappa$-complete filter, the absence of disjunctions in the language of $\kappa$-regular logic allows to derive the corresponding version of \L{}o\'s theorem by using only an appropriate $\kappa$-complete filter. The existence of this latter does not need to invoke any large cardinal assumptions. This together with the topos-theoretic characterisation of saturated models provides enough compactness to prove Grossberg conjecture.

  Although there is a vast literature covering attempts to prove Shelah's eventual categoricity conjecture (there are more than 2000 pages of approximations), our approach differs from them all in that we have a completeness theorem available, which makes possible, as we said, an omitting types theorem and its consequences. Nevertheless there are several partial results that can be found already in the literature and are derived here with fundamentally different proofs. The closest to our arguments seems to be \cite{shma}, where the strategy followed has many points in common with ours. We can point, though, the following improvements (besides the fact that we avoid large cardinals):

\begin{itemize}
\item Several known facts about AEC's, like the properties described in our Lemma \ref{lc}, Theorem \ref{cct}, Lemma \ref{at} or Theorem \ref{afsc}, are here proved in the more general context of accessible categories with directed colimits. This shows that the notion of internal size is the correct framework in which to state these results.

\item The use of the monster model is avoided altogether by using instead the syntactic category of the corresponding theory, whose universal property is enough to derive many of the statements in \cite{shma} without having to go through the semantic detour.

\item Our Theorem \ref{btwn} and its proof should correspond roughly to Proposition 4.27 of \cite{shma} (as explained in our Remark \ref{rmkesm}), but the heavy use of stability theory is here avoided through an argument relying instead on Kripke models and the completeness of infinitary intuitionistic logic, which simplifies the proof considerably.
\end{itemize}

  As we mentioned, when restricting to an AEC we recover a proof, in $ZFC$, of Shelah's eventual categoricity conjecture. Moreover, the threshold cardinal is the maximum of the Hanf numbers\footnote{The Hanf number for a certain property is the least cardinal $\mu$ such that if the property holds at a certain $\lambda \geq \mu$, then it holds in a proper class of cardinals.} for categoricity and non-categoricity, which is in fact the least possible cardinal for which the conjecture is true. We avoid $GCH$, which is spared through a forcing argument that adapts well in this particular proof, circumventing thus the fact that categoricity is not absolute. Indeed, almost all uses of $GCH$ in our paper derive from the equality $\kappa^{<\kappa}=\kappa$ for regular cardinals $\kappa$. This is in contrast with the use of it in previous literature as a device that, like uses of the diamond principle, makes the set-theoretical universe better behaved. In our case, this particular cardinal arithmetic equality will be easy to derive through the use of forcing (collapsing $\kappa^{<\kappa}$ to $\kappa$), in a way that does not modify the essential properties that we are proving (as this forcing does not change the category of models of size less than $\kappa$). This is what in the end will allow to eliminate many uses of $GCH$ in the case of AEC's, and to downgrade them to mere $SCH$ in the case of accessible categories with directed colimits. If we restrict to those accessible categories in which all morphisms are monomorphisms (in particular, if we work with AEC's), there is also an alternative argument, which we explain in the last section, that allows to achieve this same goal without recurring to forcing, by using instead the L\"owenheim-Skolem property. Thus, the reader who is not acquainted with forcing techniques can easily skip to the very last portion of the last section, in which we explain how the version of the completeness theorem that we use in our arguments does not need $GCH$  (see Remark \ref{alt}).

  Our main theorem shows that a classification theory for infinite quantifier theories is indeed possible, at least when the category of models has directed colimits. Thus, our results open up an entire new path so far unexplored in the model theory of non-elementary classes.

  The paper is structured as follows. In the section of preliminaries we introduce the notion of internal size from \cite{rlv} and establish the existence spectrum of accessible categories with directed colimits with this notion of size (which coincides in AEC's with the size of the underlying set). We then present a completely categorical account of the completeness theorem for intuitionistic logics, along with the downward L\"owenheim-Skolem theorem, which will be used for the downward categoricity transfer later. We include also the representation of a $\kappa$-classifying topos as a topos of equivariant sheaves. This latter is only needed for the proof of the generalized conjecture for accessible categories with directed colimits, but plays no r\^ole in the case of AEC's; thus the reader interested only in understanding the proof of Shelah's eventual categoricity conjecture for AEC's can disregard this material. We subsequently move to prove a technical result regarding $\kappa$-classifying toposes, introduced by the author in \cite{espindolad} as a generalization of the usual notion of classifying topos that provides the adequate context to study infinitary logics; the result established in this section is essential for the rest of the paper. The next immediate application is the omitting types theorem for infinitary theories whose categories of models have directed colimits; we provide two formulations, a syntactic one and a semantic one which does not presuppose any proof system. This theorem is in turn then used in the next section to establish a topos-theoretic characterization of $\mathcal{L}_{\infty, \kappa}$-elementary equivalence. Next, we give a precise description of the $\kappa$-classifying topos for $\kappa$-saturated models, and then complement with results on amalgamation, establishing Grossberg's conjecture in accessible categories whose morphisms are monomorphisms. Once we have amalgamation, we then present the downward categoricity transfer and the main result of the paper, namely the generalization of Shelah's eventual categoricity conjecture to accessible categories with directed colimits. We postpone for the last section the presentation of the arguments used to eliminate many of the uses of $GCH$ made thus far, which are based on a common technique described in detail in one particular case. This means that most uses of $GCH$ made before are not essential, but we prefer to prove this later to make the reading more fluid.

  To make the paper as self-contained as possible for readers not familiar with the general results of the author in \cite{espindola} and with those in the sequel \cite{espindolad}, we will include a brief description of $\kappa$-coherent theories (a particular case of $\kappa$-geometric from \cite{espindolad}), together with the notion of $\kappa$-classifying topos. In particular, we will mention the transfinite transitivity rule, which makes the completeness theorem possible. 

\section{Preliminaries}

  We start with the categorical definition of an abstract elementary class:

\begin{defs}
An abstract elementary class (AEC) $\mathcal{K}$ is a category equivalent to an accessible category with directed colimits whose morphisms are monomorphisms, that admits an embedding $F: \mathcal{K} \to \mathcal{A}$ into a finitely accessible category preserving directed colimits and monomorphisms which is, in addition:

\begin{enumerate}
\item Full on isomorphisms: for every isomorphism $h: F(A) \to F(B)$ there is an isomorphism $m: A \to B$ with $F(m)=h$.
\item Nearly full: for every commutative triangle:
\begin{center}
\begin{tikzcd}
F(A) \arrow[rd, "h"'] \arrow[rr, "F(f)"] &                          & F(B) \\
                                         & F(C) \arrow[ru, "F(g)"'] &     
\end{tikzcd}
\end{center}
there is $m: A \to C$ with $F(m)=h$.
\end{enumerate}
\end{defs}

  A result of \cite{br} gives us that an AEC automatically admits an iso-full, nearly full embedding $E: \mathcal{K} \to Emb(\Sigma)$ (for some signature $\Sigma$, where $Emb(\Sigma)$ is the category of $\Sigma$-structures and embeddings) preserving directed colimits. It also has eventually a Löwenheim-Skolem number $\lambda$: for every substructure $i: A \to F(B)$ there is $E(f): F(C) \to F(B)$ such that $i$ factors through $E(f)$ and $|E(C)| \leq |A|+\lambda$. This gives the connection with the standard model-theoretic definition of an AEC. As we remarked earlier, an important fact in AEC's is that for each object $A$, $|E(A)|$ coincides with its internal size $|A|$. This means that our eventual categoricity result for accessible categories with directed colimits is really a direct generalization of the same result for AEC's.

  When in the definition of AEC we allow for concrete $\mu$-directed colimits instead of directed, we recover the notion of $\mu$-AEC. As proven in \cite{bgrlv}, this concept admits a much simpler categorical description:

\begin{defs}
A $\mu$-abstract elementary class ($\mu$-AEC) $\mathcal{K}$ is a category equivalent to an accessible category whose morphisms are monomorphisms.
\end{defs}

  It is well-known that accessible categories are axiomatizable with sequents in $\kappa$-coherent logic (introduced in section \ref{kc}) for some $\kappa$; moreover, one can choose the signature to contain only a binary relation symbol. This changes the signature (in fact it introduces multiple sorts), but since in AEC's the size of a model (i.e. the cardinal sum of the underlying sets of each sort) coincides with its internal size, and this latter is a categorical invariant, the change of signature does not affect the categoricity spectrum; in fact, the category of models with the new signature will be equivalent to the original AEC. We will see soon that this will determine, up to equivalence, the $\kappa$-classifying topos. If one allows an infinitary signature, other axiomatizations are possible as well, as we will see later. Due to the equivalence between an accessible category and the category of models of its axiomatization, we will usually identify the object of the accessible category with the corresponding model, and loosely speak of a model and make reference to its underlying set when we deal with objects of the accessible category.

  Throughout the rest of the paper we will deal with a large accessible category $\mathcal{K}$ with directed colimits equivalent to the category of models of some $\mathcal{L}_{\kappa, \theta}$-theory. We will work with the notion of internal size $|A|$ from \cite{rlv}, defined as follows:

\begin{defs} If $r(A)$ is the least regular cardinal $\lambda$ such that $A$ is $\lambda$-presentable, then the internal size $|A|$ is defined as:

$$|A|=\left.
\begin{cases}
\kappa & \text{ if $r(A)=\kappa^+$}\\
r(A) &\text{ if $r(A)$ is limit}
\end{cases}
\right.$$
\end{defs}

We will also assume the Singular Cardinal Hypothesis ($SCH$), which states that for every infinite singular cardinal $\lambda$, we have that $\lambda^{cf(\lambda)}= 2^{cf(\lambda)} + \lambda$. This is a weakening of the Generalized Continuum Hypothesis ($GCH$), which holds, e.g., above a strongly compact cardinal, and is still powerful enough to imply well-behaved cardinal arithmetic properties. For example, it will imply (see \cite{jechst}) that $\kappa^{<\theta}=\kappa$ for all $\kappa > 2^{\theta}$ of cofinality at least $\theta$. This allows to apply the downward L\"owenheim-Skolem theorem to find submodels of size exactly $\kappa$.

  Assuming $SCH$, we can identify in an accessible category of models of some theory in $\mathcal{L}_{\kappa, \theta}$ for which cardinals the notion of internal size and that of external size (i.e., the cardinality of the underlying set of the model) coincide. Indeed, by Theorem 4.13 of \cite{rlv}, that will be the case for all regular cardinals $\lambda$ which are not successors of cardinals of cofinality less than $\theta$. This motivates the following:

\begin{defs}
$S$ is the class of cardinals which are of cofinality at least $\theta$ but are not successors of cardinals of cofinality less than $\theta$.
\end{defs}

By Theorem 4.12 of \cite{rlv} we have:

\begin{rmk}
Internal and external sizes coincide for all high enough $\lambda \in S$.
\end{rmk}

By the results of Beke and Rosicky (\cite{br}), any large accessible category with directed colimits is eventually $\lambda$-accessible for every $\lambda$. We will see in Lemma \ref{inj} below that under a categoricity assumption the subcategory where we consider only monomorphisms is closed under directed colimits above the categoricity cardinal. It follows from Corollary 4.15 of \cite{rlv} that in those conditions, and under $SCH$, there is an object of every high enough internal size whose cofinality is at least $\theta$. The following lemma shows that more can be said:

\begin{lemma}\label{singular}
Assume $SCH$ and that all morphisms are monomorphisms. If $\kappa$ is a singular cardinal of cofinality less than $\theta$, models of internal size $\kappa$ are precisely the directed colimits of models of internal size $\nu < \kappa$ in $S$.
\end{lemma}

\begin{proof}
It is easy to see that directed colimits of models of internal size $\nu < \kappa$ in $S$ must be $\kappa^+$-presentable. On the other hand, by L\"owenheim-Skolem theorem, for any subset $A$ of size $\kappa=\cup_{\nu<\kappa, \nu \in S} \nu$ of a model $M$ of internal size $\kappa$ there is a submodel of internal size $\kappa$ containing $A$. Indeed, writing $A$ as a union of a chain of $\kappa$ subsets of size less than $\kappa$, we can find submodels containing these subsets, that are of the same internal size and form moreover a chain under inclusion. Then the directed colimit of this chain will be the desired submodel. This argument also shows, as is not difficult to see, that $M$ is a $\kappa^+$-directed colimit of submodels which are themselves directed colimits of models of internal size less than $\kappa$. Since $M$ is $\kappa^+$-presentable, it follows that it must coincide with one of these submodels, as we wanted.
\end{proof}

\begin{cor}
Assume $SCH$. Then in any large accessible category with directed colimits, categorical in some cardinal, there is eventually an object of every high enough internal size. 
\end{cor}

\begin{proof}
It is enough to note that the internal size remains unchanged when one restricts to the subcategory where all morphisms are monomorphisms.
\end{proof}

In the next sections we will assume that $SCH$ holds, at least. $\mathcal{K}_{\kappa}$ will denote the subcategory of objects of internal size $\kappa$ and $\mathcal{K}_{\geq \kappa}$ that of objects of internal size at least $\kappa$.

\section{The completeness theorem}

\subsection{$\kappa$-coherent logic}\label{kc}

Let $\kappa$ be a regular cardinal such that $\kappa^{<\kappa}=\kappa$. The syntax of $\kappa$-coherent logic consists of a (well-ordered) set of sorts and a set of function and relation symbols, these latter together with the corresponding type, which is a subset with less than $\kappa$ many sorts. Therefore, we assume that our signature may contain relation and function symbols on $\gamma<\kappa$ many variables, and we suppose there is a supply of $\kappa$ many fresh variables of each sort. Terms and atomic formulas are defined as usual, and general formulas are defined inductively according to the following:

\begin{defs} If $\phi, \psi, \{\phi_{\alpha}: \alpha<\gamma\}$ (for each $\gamma<\kappa$) and $\{\psi_{\alpha}: \alpha<\delta\}$ (for each $\delta<\kappa^+$) are $\kappa$-coherent formulas, the following are also formulas: $\bigwedge_{\alpha<\gamma}\phi_{\alpha}$, $\exists_{\alpha<\gamma} x_{\alpha} \phi$ (also written $\exists \mathbf{x}_{\gamma} \phi$ if $\mathbf{x}_{\gamma}=\{x_{\alpha}: \alpha<\gamma\}$) and $\bigvee_{\alpha<\delta}\psi_{\alpha}$, this latter provided that $\cup_{\alpha<\delta}FV(\psi_{\alpha})$, the set of free variables of all $\psi_{\alpha}$, has cardinality less than $\kappa$.
\end{defs}
 
We use sequent style calculus to formulate the axioms of $\kappa$-coherent logic, as explained in \cite{espindolad}, where the system for $\kappa$-coherent logic is described. The only substantial difference with usual geometric logic (which can be seen to be exactly $\omega$-geometric logic) is the introduction of the transfinite transitivity rule (called ``rule $T$" in \cite{espindolad}): 

\begin{mathpar}
\inferrule{\phi_{f} \vdash_{\mathbf{y}_{f}} \bigvee_{g \in \gamma^{\beta+1}, g|_{\beta}=f} \exists \mathbf{x}_{g} \phi_{g} \\ \beta<\kappa, f \in \gamma^{\beta} \\\\ \phi_{f} \dashv \vdash_{\mathbf{y}_{f}} \bigwedge_{\alpha<\beta}\phi_{f|_{\alpha}} \\ \beta < \kappa, \text{ limit }\beta, f \in \gamma^{\beta}}{\phi_{\emptyset} \vdash_{\mathbf{y}_{\emptyset}} \bigvee_{f \in B}  \exists_{\beta<\delta_f}\mathbf{x}_{f|_{\beta +1}} \bigwedge_{\beta<\delta_f}\phi_{f|_{\beta+1}}}
\end{mathpar}
\\
for each cardinal $\gamma \leq \kappa$, where $\mathbf{y}_{f}$ is the canonical context of $\phi_{f}$, provided that, for every $f \in \gamma^{\beta+1}$,  $FV(\phi_{f}) = FV(\phi_{f|_{\beta}}) \cup \mathbf{x}_{f}$ and $\mathbf{x}_{f|_{\beta +1}} \cap FV(\phi_{f|_{\beta}})= \emptyset$ for any $\beta<\gamma$, as well as $FV(\phi_{f}) = \bigcup_{\alpha<\beta} FV(\phi_{f|_{\alpha}})$ for limit $\beta$. Here $B \subseteq \gamma^{< \kappa}$ consists of the minimal elements of a given bar\footnote{A bar over the tree $\gamma^{< \kappa}$ is an upward closed subset of nodes intersecting every branch of the tree.} over the tree $\gamma^{< \kappa}$, and the $\delta_f$ are the levels of the corresponding $f \in B$. 

  Basically, what it expresses is that given the tree $\gamma^{< \kappa}$ (i.e., the poset of functions $f: \beta \to \gamma$ for $\beta \leq \gamma$ with the order given by inclusion), and assuming there is an assignment of formulas $\phi_f$ to each node of the tree in such a way that the formula assigned to each node is ``covered" by the formulas assigned to its immediate successors, and the formula assigned to a node in a limit level is equivalent to the meet of the formulas assigned to its predecessors, then the formula assigned to the root must be ``covered" by the formulas assigned to the nodes ranging among the minimal elements of a given bar over the tree.  This rule is the syntactic counterpart of the corresponding exactness property of the category of sets, in which we identify each formula in a node with the set of elements where the formula holds, and where each ``cover" is really a jointly epic family of functions. 

  A $\kappa$-coherent theory corresponds precisely to a $\kappa$-coherent category, which is a coherent category with $\kappa$-small limits, stable unions of at most $\kappa$ many subobjects and satisfying that the transfinite composites of jointly epic families of arrows form a jointly epic family. Indeed, each such theory $\mathbb{T}$ gives rise to its $\kappa$-coherent syntactic category $\mathcal{C}_{\mathbb{T}}$ containing a generic model $M_0$ of $\mathbb{T}$, characterized by the universal property that models $M$ in any other $\kappa$-coherent category $\mathcal{D}$ correspond precisely to functors $\mathcal{C}_{\mathbb{T}} \to \mathcal{D}$ preserving the $\kappa$-coherent structure, i.e., to $\kappa$-coherent functors:

\begin{center}
\begin{tikzcd}
\mathcal{C}_{\mathbb{T}} \arrow[rrr, "\kappa-coherent", dashed] &  &  & \mathcal{D}          \\
M_0 \arrow[rrr, dashed, maps to] \arrow[u, maps to]          &  &  & M \arrow[u, maps to]
\end{tikzcd}
\end{center}

  Presheaves inherit the transfinite transitivity property of $\Sets$, as do sheaves over a $\kappa$-coherent category with a $\kappa$-topology (i.e., a Grothendieck topology in which transfinite composites of covers are still covers). Such toposes will be here called $\kappa$-toposes ($\kappa$-geometric toposes in \cite{espindolad}). The ($2$-)category of $\kappa$-toposes has as $1$-cells the $\kappa$-geometric morphisms, which are geometric morphisms whose inverse image preserve, in addition, all $\kappa$-small limits.

  Each $\kappa$-coherent theory admits a $\kappa$-classifying topos, introduced in \cite{espindolad} as the $\kappa$-topos containing a generic model of the theory which is universal among models in other $\kappa$-toposes, in the sense that the diagram:

\begin{center}
\begin{tikzcd}
\mathcal{C}_{\mathbb{T}} \arrow[rdd, "M"'] \arrow[rr, "Y"] &             & {\mathcal{S}et[\mathbb{T}]_{\kappa}} \arrow[ldd, "\kappa-\text{small limit preserving}", dashed] \\
                                                           & \cong       &                                                                                                  \\
                                                           & \mathcal{E} &                                                                                                 
\end{tikzcd}
\end{center}

where $\mathcal{E}$ is a $\kappa$-topos and $M$ is a $\kappa$-coherent functor from the syntactic category $\mathcal{C}_{\mathbb{T}}$, commutes up to invertible $2$-cell for an essentially unique geometric morphism (the dashed inverse image). To simplify the notation, we write throughout $(\mathcal{C}_{\mathbb{T}})_{\lambda}$, with regular $\lambda$, for the syntactic category of the theory with the axioms of $\mathbb{T}$ but in $\lambda$-coherent logic $\mathcal{L}_{\lambda^+, \lambda}$. When $\lambda$ is limit, the we understand by the notation the colimit in $\mathcal{C}at$ of $(\mathcal{C}_{\mathbb{T}})_{\kappa}$ for regular $\kappa<\lambda$. Likewise, the $\lambda$-classifying topos of $\mathbb{T}$ is written as $\Sets[\mathbb{T}]_{\lambda}$, understanding that this latter is, for limit $\lambda$, defined as the pseudolimit of the toposes $\Sets[\mathbb{T}]_{\kappa}$ for regular $\kappa<\lambda$.

  Finally, we briefly mention the $\kappa$-regular fragment as defined in \cite{makkai}, which is the subfragment of $\kappa$-coherent logic in which we remove disjunctions from the language. In this case, the transfinite transitivity rule is only applied to linear trees, where each node has just one successor, and thus reduces to the rule of dependent choices.

\subsection{Completeness of $\kappa$-coherent logic}\label{comp}

  The main tool in the analysis of infinitary logic will be the completeness theorem proved by the author in \cite{espindola} and \cite{espindolad}. Here we will present the same proof, but cast entirely in categorical language to better suit our purposes. In topos-theoretic language, it asserts that every $\kappa$-separable\footnote{A $\kappa$-separable $\kappa$-topos is one for which there is a site of size at most $\kappa$ where the $\kappa$-Grothendieck topology is generated by at most $\kappa$ many covering families.} $\kappa$-topos has enough $\kappa$-points. The $\kappa$-separability condition puts a restriction on the cardinality of formulas of the internal language of the site, which is achieved by working in the so called $\kappa$-fragments.

   In the classical case, a $\kappa$-fragment of $\mathcal{L}_{\kappa^+, \kappa}$ will be a subset of formulas formed in a language with a signature of cardinality at most $\kappa$ and a supply of $\kappa$ many fresh variables, that in addition is closed under $\kappa$-small conjunctions, disjunctions and quantification, negation and formal negation, subformulas and substitution. It follows that for any sentence $\phi$ of $\mathcal{L}_{\kappa^+, \kappa}$ there is a smallest $\kappa$-fragment containing $\phi$, and it has $\kappa$ many formulas. The same is true for any theory with at most $\kappa$ many axioms.

  An analogous definition for $\kappa$-coherent logic yields that the syntactic category of a theory with at most $\kappa$ axioms has a subcategory given by the formulas in context belonging to the $\kappa$-fragment generated by the theory. This subcategory is itself a $\kappa$-coherent category of size $\kappa$, for which we will now prove a completeness theorem.  

  Given a $\kappa$-coherent category $\mathcal{C}$ of size $\kappa$ constructed as above, we will find a jointly conservative set of $\kappa$-coherent functors $F_i: \mathcal{C} \to \Sets$. The definition of these functors is as follows. First, for each proper monomorphism $n: Y \to X$ in $\mathcal{C}$ we consider the slice $\mathcal{C}/X$. In this slice, we define for each object $B$ a set of $\kappa$ many jointly epic families of at most $\kappa$ many arrows each, in the following way. For each $A$ in $\mathcal{C}' := \mathcal{C}/X$ we consider all jointly epic families of at most $\kappa$ many arrows each generated, as a $\kappa$-Grothendieck topology, by the axioms of the theory. We then pullback these covers of $A$ along the morphism $B \to 1$, which gives us covers of $A' := A \times B$. Now we consider all sections $s: B \to A'$ of the projection $\pi_B: A' \to B$, and we pullback the $\kappa$ many covers over $A'$ along each of these sections. This defines a sequence of $\kappa$ many covers over $B$, $\mathcal{F}(B)$, which without loss of generality we can assume well-ordered and of order type $\kappa$. We then build a tree of height $\kappa$ each of whose nodes is an object of $\mathcal{C}'$ such that its immediate successors are the domains of a certain jointly epic family of arrows over the object. At step $0$, we start with the terminal object in $\mathcal{C}'$. We fix a well-ordering $w: \kappa \times \kappa \to \kappa$ such that $w(\alpha, \beta) \geq \beta$ (e.g., the canonical well-ordering of $\kappa \times \kappa$). At step $\gamma+1$ we consider $w^{-1}(\gamma)=(\alpha, \beta)$ and we pull back the $\alpha$-cover in each $\mathcal{F}(m)$, for $m$ a node in the $\beta$-level of the tree, along all the branches of height $\gamma$ over each such node $m$, to form the objects in the level $\gamma+1$ of the tree. At limit ordinals $\delta$ we take the limit of the chain of objects so defined along each branch of height $\delta$. This defines a tree of height $\kappa$ whose branches are transfinite chains of morphisms $b: 1 \leftarrow B_0 \leftarrow B_1 \leftarrow ...$. To finish the definition of the jointly conservative family of functors, we compute the corresponding slices $\mathcal{C}' \to \mathcal{C}'/B_0 \to \mathcal{C}'/B_1 \to ...$ with (some choice of) pullback functors between them, and take their (pseudo-)colimit in $\mathcal{C}at$, $\mathcal{C}'_b$. Our family of functors will then be the composites $\mathcal{C} \to \mathcal{C}' \to \mathcal{C}'_b \to \Sets$, where we vary over all proper monomorphisms $n: Y \to X$ and all branches $b$, and where each $\mathcal{C}'_b \to \Sets$ is given by the term model functor $[1, -]$. 

  We need to prove that each functor is $\kappa$-coherent, for which it is enough to prove that each one of these composites $[1, -] \circ F$ preserves jointly epic families of at most $\kappa$ many arrows. Suppose $\{C_i \twoheadrightarrow M_i \rightarrowtail A\}_{i<\kappa}$ is such a family and we have in $\mathcal{C}'_b$ an arrow $s: 1 \to \bigvee_{i<\kappa}F(M_i) =F(A)$, which is a section of the map $F(A) \to 1$. By construction, such an arrow has a representative $s_B: B \to A'$ in some slice $\mathcal{C}'/B$ which is a section of a map $h: A' \to B$, so that, by definition, the image $\{C_i' \twoheadrightarrow M_i' \rightarrowtail A'\}_{i<\kappa}$ in $\mathcal{C}'/B$ of the family $\{C_i \twoheadrightarrow M_i \rightarrowtail A\}_{i<\kappa}$ , when pulled back along the section $s_B: B \to A'$, provides a jointly epic family $\{P_i \to B\}_{i<\kappa}$ that belongs to $\mathcal{F}(B)$. Such a jointly epic family is then, in turn, pulled back at some step in the transfinite construction to form the successors of some node $C$, so that the branch $b$ chosen for our functor $F$ will have some object $P_i' := P_i \times_B C$ as a successor of $C$. It follows that in the slice $\mathcal{C}'/P_i'$ the section $s_{P_i'}: P_i' \to Q$, which represents $s$, will factor through $i: P \to Q$ via the induced morphism $(f, Id_{P_i'}): P_i' \to P$ to the pullback $P$ (which is the image of $C_i$ in the slice $\mathcal{C}'/P_i'$). That this is the case can be seen by noticing that $s_{P_i'}: P_i' \to Q$ and $i \circ (f, Id_{P_i'}): P_i' \to Q$ give the same morphism when composed with $a$ and $b$ (see the diagram below).

\begin{center}
\begin{tikzcd}
               &                           &                              &                           &                                                         & P \arrow[ld] \arrow[rr]         &                                                  & Q \arrow[ld, "a"'] \arrow[dd, "b", bend right]                                                   \\
               & C_i' \arrow[ld] \arrow[r] & A' \arrow[ld] \arrow[d, "h"] &                           & C' \times_{B}C \arrow[ld] \arrow[rr]                    &                                 & A' \times_{B}C \arrow[ld] \arrow[dd, bend right] &                                                                                                  \\
C_i \arrow[r] & A \arrow[d]               & B \arrow[ld]                 & C_i' \arrow[rr]           &                                                         & A' \arrow[dd, "h"', bend right] &                                                  & P_i' \arrow[ld, "g"] \arrow[uu] \arrow[lluu, "{(f, Id_{P_i'})}", dashed, bend right, shift left] \\
               & 1                         &                              &                           & P_i' \arrow[rr, "\qquad g"] \arrow[uu, "f "] \arrow[ld] &                                 & C \arrow[ld] \arrow[uu]                          &                                                                                                  \\
               &                           &                              & P_i \arrow[rr] \arrow[uu] &                                                         & B \arrow[uu]                    &                                                  &                                                                                                 
\end{tikzcd}
\end{center}

  Therefore, in $\mathcal{C}'_b$, $s: 1 \to F(\bigvee_{i<\kappa}M_i) = F(A)$ factors through some $F(C_i)$, as we wanted to prove.

  Finally, it should be clear now that this family of $\kappa$-coherent functors is jointly conservative, since for every proper monomorphism $m: Y \to X$, its image in the slice $\mathcal{C}' := \mathcal{C}/X$ is a proper subterminal object. If all functors in our family sent this to an isomorphism, it would be an isomorphism in each $\mathcal{C}'_b$ and hence also in some slice $\mathcal{C}'/B_b$ for some $B_b$ in each branch $b$. Now the transfinite transitivity property implies that the family of functors $\{\mathcal{C}' \to \mathcal{C}'/B_b\}_{b}$ is jointly conservative, getting that the original monomorphism was not proper after all. Whence, the family really is jointly conservative, as we wanted.

  Note that the cardinal arithmetic assumption on $\kappa$ that we need to put is $\kappa^{<\kappa}=\kappa$ (in particular, $\kappa$ must be regular). This guarantees that the set of $\kappa$-coherent formulas in a $\kappa$-fragment containing any set of $\kappa$ many of them will have cardinality $\kappa$, so that each set of covers $\mathcal{F}(B)$ has size $\kappa$ and we can perform the transfinite construction up to $\kappa$. There is however the possibility to extend the theorem even if $\kappa$ is singular, provided the fragment in which the axioms of our theory are expressed is the union of all $\nu$-fragments of size $\nu$ for regular $\nu<\kappa$ and that we restrict to those morphisms which are monomorphisms. This happens because then the well-ordering of $\kappa \times \kappa$ can be chosen to be the union of the canonical well-orderings of each $\nu \times \nu$. This observation will be used later on when we prove the downward categoricity transfer.

\subsection{The L\"owenheim-Skolem theorem}\label{lst}

Recall that the universal property of the slice category (which has $2$-categorical aspects) says that the pullback functor to the slice is the universal morphism into a $\kappa$-coherent category which has a section (the diagonal $\Delta:(\mathbf{x}, \phi) \to (\mathbf{x}, \phi) \times (\mathbf{x}, \phi)$) for the morphism $(\mathbf{x}, \phi) \to 1$):

\begin{center}
\begin{tikzcd}
\mathcal{C}_{\mathbb{T}} \arrow[rr] \arrow[rdd, "F"'] &                           & {\mathcal{C}_{\mathbb{T}}/(\mathbf{x}, \phi), \Delta} \arrow[ldd, "\overline{F} \qquad \overline{F}(\Delta)=\mathbf{c}", dashed] \\
                                                      & \cong                     &                                                                                                                                  \\
                                                      & {\mathcal{D}, \mathbf{c}} &                                                                                                                                 
\end{tikzcd}
\end{center}

This universal property allows us to have an entirely categorical understanding of the completeness theorem just proven, since starting with the syntactic category of a $\kappa$-coherent theory $\mathbb{T}$, each slice over $(\mathbf{x}, \phi)$ is equivalent to the syntactic category of the theory $\mathbb{T} \cup \phi(\mathbf{c})$ where $\mathbf{c}$ is a tuple of fresh constants. 

 Since the construction explained in the proof of the completeness theorem proceeds by successively taking slices over the syntactic category $\mathcal{C}$ of $\mathbb{T}$, an easy application of the universal property shows that for any $\kappa$-coherent model $M: \mathcal{C} \to \Sets$ and any tuple of less than $\kappa$ elements $\mathbf{c} \subseteq M$ (where we identify the functor with the underlying set of $M((x, \top))$), factors, up to invertible $2$-cell, through some $\mathcal{C}'_b$ for some branch $b$. Indeed, we know it factors through the slice $\mathcal{C}' := \mathcal{C}/(\mathbf{x}, \top)$; since $M$ is $\kappa$-coherent, it sends the tree of height $\kappa$ built by transfinite recursion to a tree of sets each of whose nodes has a jointly epic family of morphisms from its immediate successors. This allows to define successively a compatible chain of elements $\mathbf{d}_B \in M(B)$ which, by the universal property of the slice, induce morphisms $M_B: \mathcal{C}'/B \to \Sets$ along objects $B$ in some branch $b$, with invertible $2$-cells between them that form a pseudococone. It follows, as we claimed, that then there is an induced $\kappa$-coherent functor $\overline{M}: \mathcal{C}'_b \to \Sets$ that extends $M$. Since $M$ must be a colimit of representable functors in whose diagram there is always $[1, -]$, this latter representable provides, as we have proven above, a $\kappa$-coherent submodel $M'$ of $M$ containing the tuple $\mathbf{c}$. In other words, we have proven a version of the downward L\"owenheim-Skolem theorem which will become essential for the downward categoricity transfer, as we will see.

  Note that the submodel $[1, -] \circ F$ is a $\kappa$-coherent model of each theory $\mathbb{T} \cup \phi(\mathbf{c})$, where $F((\mathbf{x}, \phi)) = B_i$ is the $i$-th object in the branch $b$. This readily implies that it must also be a model of the subcategory of the syntactic category of the $\kappa$-coherent theory axiomatized by $\mathbb{T} \cup \{\phi_i(\mathbf{c_i})\}_{i<\kappa}$ which contains formulas in the $\kappa$-fragment generated by these axioms. In turn, such a subcategory has an induced canonical morphism from $\mathcal{C}'_b$ due to the universal property of this latter. This will be used for the Kripke completeness theorem.

\subsection{Kripke completeness for $\kappa$-first-order theories}\label{iil}

  The categorical proof of the completeness theorem can be applied in particular to any $\kappa$-Heyting category, e.g. the syntactic category of a $\kappa$-Heyting theory $\mathbb{T}$ of size $\kappa$. The $\kappa$-Heyting structure is preserved under slicing; moreover, the colimit $\mathcal{C}'_b$ has a functor $[1, -]$ such that the composite $[1, -] \circ F$ is a model of the subcategory of the $\kappa$-Heyting syntactic category of the theory $\mathbb{T} \cup \{\phi_i(\mathbf{c_i})\}_{i<\kappa}$ which contains formulas in the $\kappa$-fragment generated by these axioms (see the last paragraph of section \ref{lst}). In turn, the syntactic category of this latter theory has, by the universal property of $\mathcal{C}'_b$, an induced morphism from it. Hence, the construction really gives us a set of $\kappa$-Heyting prime theories\footnote{A prime theory is a deductively closed theory with the disjunction and existence properties.} over the language in which we add $\kappa$ many fresh constants (namely, the $\kappa$-Heyting theories $\mathbb{T}_b := \mathbb{T} \cup \{\phi_i(\mathbf{c_i})\}_{i<\kappa}$). These have, in addition, a conservativity property: given any tuple of less than $\kappa$ constants $\mathbf{c}$ in the new language such that $\mathbb{T}, \phi(\mathbf{c}) \nvdash \psi(\mathbf{c})$, there is one theory $\mathbb{T}_b$ in our set of new $\kappa$-Heyting theories such that $\mathbb{T}_b, \phi(\mathbf{c}) \nvdash \psi(\mathbf{c})$. As it is well-known, at least in the finitary case, this property is essentially what allows to build a set of jointly conservative Kripke models\footnote{A Kripke model of $\mathbb{T}$ over $P$ is a $\kappa$-Heyting functor $\mathcal{C}_{\mathbb{T}} \to \Sets^P$ for a given poset $P$.} of $\mathbb{T}$ over a tree of height $\omega$. Indeed, we can iterate the construction that builds the set of prime theories over the extended language $\omega$ many times, continuing with each $\mathcal{C}_{\mathbb{T}_b}$ at limit levels, forming the said tree. Each node is identified with the term model of the corresponding $\mathcal{C}'_b$, with the induced homomorphisms between them arising from the L\"owenheim-Skolem theorem (cf. section \ref{lst}). It is not hard to see that the theory $\mathbb{T}_b$, in the appropriate language, represents the set of all $\kappa$-Heyting sentences that are forced in the corresponding node. To prove this, we can make a straightforward induction on the complexity of the formula using the property of the set of prime theories (extending the theory of each node) that we have mentioned above. Since the functors between each category $\mathcal{C}_{\mathbb{T}_b}$ and each of its extensions is $\kappa$-Heyting, a formula proved by $\mathbb{T}_b$ must be forced in the corresponding node. On the other hand, the theories $\mathbb{T}_b$ are prime theories for the following reason: if $\mathbb{T}_b$ proves $\bigvee_i\exists \mathbf{x}\phi_i$, its $\kappa$-coherent Morleyization\footnote{The $\kappa$-coherent Morleyization of $\mathbb{T}_b$ is the theory of $\kappa$-coherent models of $\mathcal{C}_{\mathbb{T}_b}$. Each formula $\phi$ in the original signature is then equivalent, in all $\kappa$-coherent models of $\mathcal{C}_{\mathbb{T}_b}$, to a formula $P_{\phi}$ in the new signature (its $\kappa$-coherent Morleyization).} $\bigvee_i\exists \mathbf{x}P_{\phi_i}$ holds in the $\kappa$-coherent model corresponding to its node, whence, some $P_{\phi_i}(\mathbf{c})$ must hold there and therefore also holds in the successor nodes. Thus, by the construction of section \ref{comp}, it is provable in $\mathcal{C}_{\mathbb{T}_b}$. To see that forced formulas are provable in the corresponding theory, the non trivial case is when $\phi$ is of the form $\forall \mathbf{x} (\psi \to \eta)$. Then if a node forces $\phi$, the sequent $\psi \vdash_{\mathbf{x}} \eta$ must be forced in all successor nodes, being thus proved in each successor theory by inductive hypothesis. Therefore, its $\kappa$-coherent Morleyization holds in all successor nodes, so that by the construction of section \ref{comp}, it is provable in $\mathcal{C}_{\mathbb{T}_b}$. 

 If we start with the $\kappa$-Heyting theory $\mathbb{T}$, this procedure provides a tree of prime theories over each immediate successor which will be a jointly conservative set of Kripke models for $\mathbb{T}$. We will make use of this observation for the downward categoricity transfer. 

\subsection{Representation theorem for $\kappa$-toposes with enough $\kappa$-points}\label{iaf}

  Awodey and Forssell provided in \cite{af} a logical approach to the representation theorem of Butz and Moerdijk for toposes with enough points, providing actually a reconstruction result that allows to recover the theory (up to pretopos completion) from its category of models. Unlike previous reconstruction results, which relied on using the ultraproduct structure in the category of models, this result puts instead a topological structure that allows to recover the classifying topos of the theory as the topos of equivariant sheaves on a topological groupoid. This is possible because the category of models carries a natural topology which can be presented in terms of the formulas of theory, to which one can associate the basic opens. 

  The fact that the set of opens satisfy the axioms of a topology is then a consequence of the formula construction, and in particular, the fact that a finite intersection of opens is open is directly linked to the fact that a finite conjunction of coherent formulas is coherent. When we work in infinitary logic like, e.g., $\mathcal{L}_{\kappa^+, \kappa}$, the formula construction allows for infinitary conjunctions, which in turn gives the associated topology on the category of models a distinctive property: the intersection of less than $\kappa$ many open sets is still open. We refer to it as a $\kappa$-topology. It turns out that many of the properties and results on topological spaces have also variants for $\kappa$-topologies. As a result, almost all the setup available in \cite{af} is readily generalizable to the infinitary case as soon as we apply our completeness theorem relating the syntax and the semantics. Such a completeness theorem is crucial since the classifying topos is obtained through the representation theorem of Butz and Moerdijk, which requires the topos to have enough points. 

  A similar representation theorem for $\kappa$-classifying toposes can now be shown. In the case of $\kappa$-Grothendieck toposes, it turns out that their main characteristic (the transfinite transitivity property) is shared by toposes of equivariant sheaves on topological groupoids where we use $\kappa$-topologies. This is what allows, in essence, the transfer of Butz-Moerdijk result to the infinitary case; we briefly mention the steps to do it:

\begin{enumerate}

\item The logical topology on the set of models corresponding to a $\kappa$-coherent decidable theory is a $\kappa$-topology, meaning that intersection of less than $\kappa$ open sets is open. This is evident from the fact that intersections of basic opens corresponding to $\kappa$-coherent formulas corresponds in turn to the conjunction of the formulas, which is $\kappa$-coherent. 

\item Sheaves on a $\kappa$-topological space $X$ (étale bundles over $X$) satisfy the transfinite transitivity property. This is a direct consequence of a result proven in \cite{espindolad}, namely, that when the underlying category of the site satisfies the transfinite transitivity property (as is the case with the lattice of open sets of a $\kappa$-topology), so does the sheaf topos on it.

\item Continuous functions between $\kappa$-topological spaces give rise to $\kappa$-geometric morphisms between the corresponding sheaf toposes, i.e., geometric morphisms whose inverse images preserve $\kappa$-limits. Moreover, morphisms of sites that satisfy the transfinite transitivity property induce $\kappa$-geometric morphisms. This is essentially contained in \cite{espindolad}. With a proof similar to that of the universal property of a $\kappa$-classifying topos we can see that morphisms of sites satisfying the transfinite transitivity property give rise to $\kappa$-geometric morphisms between their corresponding $\kappa$-Grothendieck toposes on them.

\item Given a $\kappa$-topological groupoid, the topos of equivariant sheaves, defined as a (pseudo-) colimit of toposes, satisfies the transfinite transitivity property. This is a consequence of the construction of pseudo-colimits of $\kappa$-Grothendieck toposes, as explained, e.g., in \cite{moerdijk2} or \cite{kelly}: it suffices to prove that the transfinite transitivity property is preserved during the construction of coproducts, tensor products, iso-coinserters and iso-coequifiers. Alternatively, we can use the explicit description of the topos of equivariant sheaves in terms of local homeomorphisms. 

\item The theory is recovered up to $\kappa$-pretopos completion as the full subcategory of the $\kappa^+$-compact decidable subobjects. (A $\kappa^+$-compact object is one such that every family of subobjects that cover it contains a subfamily of at most $\kappa$ subobjects that still covers). This is a consequence of the universal property of the $\kappa$-classifying topos.

\end{enumerate}

More explicitely: starting with the category of models of a $\kappa$-coherent theory, we define the following topological groupoid of models. Fixing a set $U$ of size $\kappa$, we consider the set $G_0$ of all models whose underlying structure has elements from $U$. Then we consider the set $G_1$ of isomorphisms between them. We put a topology on $G_0$ whose basic opens consist of sets of models of the form:

$$(\phi(\mathbf{x}), \mathbf{c}):=\{M \in G_0: \mathbf{c} \in [[\phi(\mathbf{x})]]^M\}$$
\\
Next, we consider the least topology on $G_1$ which makes the source and target $s, t: G_1 \to G_0$ continuous and that contains all sets of the form:

$$(\mathbf{a} \mapsto \mathbf{b}):=\{f \in G_1: \mathbf{a} \in s(f) \text{ and } f(\mathbf{a})=\mathbf{b}\}$$

  In our case, the topos of equivariant sheaves on our topological groupoid has as objects pairs $(a: A \to G_0, \alpha)$, where $a$ is a local homeomorphism and $\alpha: G_1 \times_{G_0} A \to A$ satisfies the conditions:
\begin{itemize}
\item $a(\alpha(f, x))=t(f)$
\item $\alpha(1, x)=x$
\item $\alpha(g, \alpha(f, x))=\alpha(gf, x)$
\end{itemize}

  We get in the end:

\begin{thm}
The $\kappa$-classifying topos of $\mathbb{T}$ is precisely the topos of equivariant sheaves $\mathcal{S}h_{G_1}(G_0)$
\end{thm}

The proof follows essentially the same outline of Theorem 1.4.8 in \cite{af}, with the only modification that we consider syntactic categories of $\kappa$-theories and topological groupoids of models with respect to the $\kappa$-logical topology. The only detail which needs a different justification is in showing that the embedding $\mathcal{C}_{\mathbb{T}} \to \mathcal{S}h_{G_1}(G_0)$ given by sending a formula to its corresponding definable set functor is cover-reflecting when we consider the canonical coverage in $\mathcal{S}h_{G_1}(G_0)$ and the $\kappa$-coherent coverage in $\mathcal{C}_{\mathbb{T}}$. But this is exactly given by the completeness theorem of section \ref{comp}, since the models of $G_0$ are enough.

\section{The $\lambda$-classifying topos of a $\kappa$-theory}

In this section fix $\kappa<\lambda$ such that $\kappa^{<\kappa}=\kappa$ and $\lambda^{<\lambda}=\lambda$. Let $\mathbb{T}$ be a $\kappa$-coherent theory in $\mathcal{L}_{\kappa^+, \kappa}$, $\mathcal{C}_{\mathbb{T}}$ be its syntactic category and $Mod_{\lambda}(\mathbb{T})$ be the full subcategory of $\lambda$-presentable models. Assume that the category of models of $\mathbb{T}$ is $\lambda$-accessible (this is the case, e.g., if $\lambda=\kappa^+$ or, more generally, if $\kappa^+ \trianglelefteq \lambda$). Let $\mathbb{T}'$ be the theory in $\mathcal{L}_{\lambda^+, \lambda}$ with the same axioms as those of $\mathbb{T}$. An important result we will prove here is the following:

\begin{thm}\label{lk}
The $\lambda$-classifying topos of $\mathbb{T}'$ is equivalent to the presheaf topos $\Sets^{Mod_{\lambda}(\mathbb{T})}$. Moreover, the canonical embedding of the syntactic category $\mathcal{C}_{\mathbb{T}'} \hookrightarrow \Sets^{Mod_{\lambda}(\mathbb{T})}$ is given by the evaluation functor, which on objects acts by sending $(\mathbf{x}, \phi)$ to the functor $\{M \mapsto [[\phi]]^M\}$.
\end{thm}

\begin{proof}
By hypothesis every model of $\mathbb{T}'$ is a $\lambda$-filtered colimit of models in $Mod_{\lambda}(\mathbb{T})$. Note first that the following diagram:

\begin{displaymath}
\xymatrix{
\mathcal{C}_{\mathbb{T}'} \ar@{^{(}->}[rr]^{ev} \ar@{->}[ddr]_{M \cong \varinjlim_i M_i} & & \Sets^{Mod_{\lambda}(\mathbb{T})} \ar@{->}[ddl]^{M' \cong \varinjlim_i ev_{M_i}}\\
 & & \\
 & \Sets & \\
}
\end{displaymath}
\noindent
commutes up to invertible $2$-cell. Here $ev$ and $ev_{M_i}$ are the evaluation functors, defined on objects as $ev((\mathbf{x}, \phi))=\{M \mapsto [[\phi]]^M\}$ and $ev_{M_i}(F)=F(M_i)$, respectively, while $\varinjlim M_i$ is the canonical $\lambda$-filtered colimit of $\lambda$-presentable models associated to the model $M$. Note also that since $\lambda$-filtered colimits commute with $\lambda$-small limits, $M'$ will preserve, in addition to all colimits, also $\lambda$-small limits.

  Let now $\Sets[\mathbb{T}']_{\lambda}$ be the $\lambda$-classifying topos of $\mathbb{T}'$. We shall prove that this latter is equivalent to $\mathcal{S}et^{Mod_{\lambda}(\mathbb{T})}$ by verifying in this presheaf topos the universal property of $\Sets[\mathbb{T}']_{\lambda}$, i.e., that models of $\mathbb{T}'$ in a $\lambda$-topos $\mathcal{E}$ corresponds to $\lambda$-geometric morphisms from $\mathcal{E}$ to the presheaf topos. It is enough to prove this universal property in the particular case in which $\mathcal{E}=\Sets[\mathbb{T}']_{\lambda}$. 

  Given then the $\lambda$-classifying topos $\mathcal{E}$ of $\mathbb{T}'$, by the completeness theorem (see section \ref{comp}) it will have enough $\lambda$-points. Hence, there is a conservative $\lambda$-geometric morphism with inverse image $E: \mathcal{E} \to \Sets^{I}$ such that composition with the evaluation at $i \in I$, $ev(i)E$ gives a $\lambda$-point of $\mathcal{E}$. Now each model of $\mathbb{T}'$ in $\mathcal{E}$, $N: \mathcal{C}_{\mathbb{T}'} \to \mathcal{E}$ gives rise to models in $\Sets$ by considering their images through each $ev(i)E$. These correspond to unique (up to isomorphism) $\lambda$-geometric morphisms with inverse image $\Sets^{Mod_{\lambda}(\mathbb{T})} \to \Sets$, which in turn induce a $\lambda$-geometric morphism with inverse image $G: \Sets^{Mod_{\lambda}(\mathbb{T})} \to \Sets^I$ and with the property that the composition $G \circ ev: \mathcal{C}_{\mathbb{T}'} \to \Sets^{Mod_{\lambda}(\mathbb{T})} \to \Sets^I$ is the same (up to isomorphism) as $EN: \mathcal{C}_{\mathbb{T}'} \to \Sets^I$. In other words, considering $\mathcal{E}$ as a subcategory of $\Sets^I$, the image of $G \circ ev$ belongs to $\mathcal{E}$.

\begin{displaymath}
\xymatrix{
\mathcal{C}_{\mathbb{T}'} \ar@{^{(}->}[rr]^{ev} \ar@{->}[ddr]_{N} & & \Sets^{Mod_{\lambda}(\mathbb{T})} \ar@{-->}[ddl] \ar@{->}[dddl]^{G} \\
 & & \\
 & \mathcal{E} \ar@{->}[d]_{E} & \\
 & \Sets^I \ar@{->}[d]_{ev(i)} & \\
  & \Sets & \\
}
\end{displaymath}
\noindent

  On the other hand, every object $F$ in $\Sets^{Mod_{\lambda}(\mathbb{T})}$ can be canonically expressed as a colimit of representables, $F \cong \varinjlim_i [M_i, -]$. In turn, and supposing for a moment that each $|M_i|<\lambda$, we have that $M_i: \mathcal{C}_{\mathbb{T}} \to \Sets$ is a $\lambda$-small colimit of representables $M_i \cong \varinjlim_j [\phi_{ij}, -]$.

It follows that:

$$F \cong \varinjlim_i[\varinjlim_j[\phi_{ij}, -]_{\mathcal{C}_{\mathbb{T}}}, -]_{Mod_{\lambda}(\mathbb{T})} \cong \varinjlim_i \varprojlim_j [[\phi_{ij}, -]_{\mathcal{C}_{\mathbb{T}}}, -]_{Mod_{\lambda}(\mathbb{T})} \cong \varinjlim_i \varprojlim_j ev(\phi_{ij})$$
\noindent
where the last isomorphism follows from Yoneda lemma. We claim that this even happens if $|M_i|=\lambda$. Indeed, in this case $\lambda$ is the successor of a singular cardinal $\mu$ of cofinality less than $\theta$ and $M_i$ has internal size $\mu$, so that by Lemma \ref{singular} it is a directed colimit of models of cardinality less than $\mu$. Suppose first $\kappa<\mu$. Then we can express:

$$[M_i, -] \cong [\varinjlim_k \varinjlim_j[\phi_{ijk}, -]_{\mathcal{C}_{\mathbb{T}}}, -]_{Mod_{\lambda}(\mathbb{T})} \cong \varprojlim_k \varprojlim_j ev(\phi_{ijk})$$
\noindent
and $F$ has a similar expression, where limits are $\lambda$-small. If, on the other hand, $\kappa=\mu$, then we note that the restriction functor $r: \Sets^{\mathcal{K}_{\geq \nu, \leq \mu}} \to \Sets^{\mathcal{K}_{\mu}}$ has a right adjoint $i$ such that $ri(F) \cong F$ for each $F$ in $\Sets^{\mathcal{K}_{\mu}}$, and expressing now $i(F)$ as a colimit of limits of evaluations, as above, we can now apply $r$, which preserves colimits and limits, and arrive to the same expression for $F$.

Now $G$ preserves $\lambda$-small limits and colimits, and so we will have:

$$G(F) \cong \varinjlim_i \varprojlim_j G \circ ev(\phi_{ij})$$
\noindent
and similarly on arrows. Therefore, $G$ is completely determined (up to isomorphism) by its value on the objects $ev(\phi_{ij})$. Since the value of $G$ on such objects belongs to $\mathcal{E}$, and $E$ preserves $\lambda$-small limits and colimits, it follows that $G$ itself factors through $\mathcal{E}$. Moreover, it is the unique (up to isomorphism) inverse image of a $\lambda$-geometric morphism corresponding to the given model in $\mathcal{E}$. This finishes the proof.
\end{proof}

\section{The omitting types theorem for infinite quantifier languages}

We will assume in this section that $\kappa^{<\kappa}=\kappa$. By a type we understand a consistent set of formulas in a given tuple of variables. It is complete when the set is maximal. The goal of this section is to prove the following:

\begin{thm}\label{ott}
(Omitting types theorem for infinite quantifier languages) Assume $GCH$, and let $\kappa$ be a regular cardinal. Let $F$ be a Boolean $\kappa$-fragment of $\mathcal{L}_{\kappa^+, \kappa}$ containing a consistent theory $\mathbb{T}$ whose category of models has directed colimits, and let $\{p_i: i < \kappa\}$ be a set of non-isolated types. Then there is a model of $\mathbb{T}$ that simultaneously omits all the types.
\end{thm}

Note that this version of the omitting types theorem can be expressed in an entirely semantical way:

\begin{thm}
Assume $GCH$. Let $\mathbb{T}$ be a satisfiable theory in a Boolean $\kappa$-fragment of $\mathcal{L}_{\kappa^+, \kappa}$ whose category of models has directed colimits, and let $p_i$, for each $i < \kappa$, be a set of formulas of the fragment. Suppose that whenever $\psi$ is such that $\mathbb{T} \cup \exists \mathbf{x} \psi$ is satisfiable, there is $\phi$ in $p_i$ such that $\mathbb{T} \cup \exists \mathbf{x}(\psi \wedge \neg \phi)$ is satisfiable. Then the theory:

$$\mathbb{T} \cup \bigwedge_{i<\kappa} \forall \mathbf{x} \bigvee_{\phi \in p_i} \neg \phi(\mathbf{x})$$

is satisfiable.
\end{thm}

\begin{proof}
Consider the (Boolean) syntactic category of $\mathbb{T}$ in $\mathcal{L}_{\kappa^+, \kappa}$ and the subcategory $\mathcal{C}_{\mathbb{T}}$ given by those formulas in context $[\mathbf{x}, \phi]$ belonging to the $\kappa$-fragment $F$. For each type $p_i=\{\phi_i(\mathbf{x})\}_{i<\kappa}$ consider the family of arrows $U_i=\{[\mathbf{x}, \neg \phi_i] \to [\mathbf{x}, \top]\}_{i<\kappa}$. Put a $\kappa$-Grothendieck topology $\tau$ on $\mathcal{C}_{\mathbb{T}}$ generated by:

\begin{enumerate}

\item all $\kappa$-small jointly epic families of arrows and the $\kappa^+$-small jointly epic families of arrows corresponding to axioms of (the $\kappa$-coherent Morleyization of) $\mathbb{T}$
\item the families $U_i$ for each type $p_i$

\end{enumerate}

It follows that a $\kappa$-flat continuous functor $\mathcal{C}_{\mathbb{T}} \to \mathcal{S}et$, i.e., a $\kappa$-point of the corresponding $\kappa$-topos of sheaves, is exactly a model of $\mathbb{T}$ omitting all of the $p_i$. This topos is clearly $\kappa$-separable, so that by the completeness theorem (see section \ref{comp}) it will have enough $\kappa$-points. However, we need to verify that it is non-degenerate to guarantee that there will be at least one non-trivial such model (it is easy to see that, if one of the types is isolated, the topos is degenerate, but we will see that this is the only obstruction).

It is enough to verify that the representable functor $[-, 0]$ is a sheaf for any $\tau$-covering family, since then the conservativity of Yoneda embedding will imply that $\mathcal{S}h(\mathcal{C}_{\mathbb{T}}, \tau)$ is non-degenerate. Now any such $\tau$-covering family is built via pullbacks and transfinite composites from the two types of covers specified above. Clearly, $[-, 0]$ is a sheaf for the first type of covers. The covering $U_i$ becomes, in the $\kappa$-classifying topos $\mathcal{E}$ of $\mathbb{T}$ (that is, when sheafifying with respect to only the first type of covers) a family $\{\neg C_i \to A\}$, not necessarily epimorphic. But since the type $p_i$ was non-isolated, it follows that $\bigwedge_{i<\kappa}C_i=0$ in $\mathcal{E}$, or, what is the same, $A=\neg \neg \bigvee_{i<\kappa} \neg C_i$ there. This means that the family $\{\neg C_i \to A\}$, while not necessarily covering, is covering up to a double negation. It is easy to see that this is also true for the image of the cover through the conservative evaluation functor $\mathcal{E} \to \Sets^{Mod(\mathbb{T}_{\kappa})}$, by using that in $\Sets^{Mod(\mathbb{T}_{\kappa})} \cong \Sets[\mathbb{T}_{\kappa}]_{\kappa^+}$ the sequent $\bigwedge_{i<\kappa}\neg \neg \phi_i \vdash_{\mathbf{x}} \neg \neg \bigwedge_{i<\kappa}\phi_i$ holds (which can be easily verified to be a consequence of the category of models having directed colimits). Indeed, this is equivalent to saying that if the pullback of the cover along a morphism $A \to (\mathbf{x}, \top)$ is $0$, then $A \cong 0$. Since conjunctions of size $\kappa$ of atomic formulas generate $\Sets^{Mod(\mathbb{T}_{\kappa})}$, it is enough to consider those $A$ of that form. Since the statement holds when $A$ is an atomic formula, it must hold still for such a conjunction due to the validity of the sequent $\bigwedge_{i<\kappa}\neg \neg \phi_i \vdash_{\mathbf{x}} \neg \neg \bigwedge_{i<\kappa}\phi_i$ and the construction of pullbacks in the syntactic category $(\mathcal{C}_{\mathbb{T}_{\kappa}})_{\kappa^+}$. In summary, our dense covers in $\mathcal{E}$ remain dense in $\Sets[\mathbb{T}_{\kappa}]_{\kappa^+}$. The same is, of course, true for pullbacks of such families, and in fact for a transfinite composite of such families (for this latter fact we use again the sequent $\bigwedge_{i<\kappa}\neg \neg \phi_i \vdash_{\mathbf{x}} \neg \neg \bigwedge_{i<\kappa}\phi_i$ and the transfinite transitivity rule from \cite{espindolad}). In particular, this means that if the domains of the arrows in a $\tau$-covering family are $0$, so is the common codomain. This says precisely that $[-, 0]$ is a sheaf for the $\tau$-covering family, as we wanted. 
\end{proof}

\section{Partial isomorphisms and the $\lambda$-classifying topos}

We will prove now a connection between $\mathcal{L}_{\infty, \lambda}$-equivalence and $\lambda$-classifying toposes which will be useful. It is essentially a consequence of the omitting types theorem we proved before (Theorem \ref{ott}), and is inspired by model-theoretic arguments of Vaught on atomic and prime models and topos-theoretic results from Blass and \v{S}\v{c}edrov on Boolean classifying toposes.

\begin{thm}\label{cct}
Let $\kappa$ be a regular cardinal such that $\kappa^{<\kappa}=\kappa$. Let $\mathbb{T}$ be a theory in a Boolean $\kappa$-fragment of $\mathcal{L}_{\kappa^+, \kappa}$ whose category of models has directed colimits. Then for any $\lambda \geq \kappa$ such that $\lambda^{<\lambda}=\lambda$, every pair of models of $\mathbb{T}$ of size $\lambda$ are $\mathcal{L}_{\infty, \lambda}$-elementarily equivalent if and only if the $\lambda$-classifying topos of the theory $\mathbb{T}_{\lambda} := \mathbb{T} \cup \{``\text{there are $\lambda$ distinct elements"}\}$ is two-valued and Boolean (alternatively, atomic and connected). 
\end{thm}

\begin{proof}
($\implies$) Suppose any two models of $\mathbb{T}$ of size $\lambda$ are $\mathcal{L}_{\infty, \lambda}$-elementarily equivalent and consider the syntactic category $\mathcal{C}$ of the theory $\mathbb{T} \cup \{``\text{there are $\lambda$ distinct elements"}\}$, axiomatized in an appropriate Boolean fragment of $\mathcal{L}_{\lambda^+, \lambda}$\footnote{This can be done through sentences expressing that for each tuple of less than $\lambda$ distinct elements there is one element different from all those of the tuple.}. This latter theory is clearly complete, since by (downward) L\"owenheim-Skolem theorem it follows that all models are $\mathcal{L}_{\lambda^+, \lambda}$-elementarily equivalent to a model $\mathcal{M}$ of cardinality $\lambda$. Therefore, its $\lambda$-classifying topos must be two-valued. To see that it is Boolean, we will prove first that $\mathcal{C}$ is atomic, i.e., each Boolean algebra of subobjects of a given object is atomic. 

Let $[\mathbf{x}, \psi(\mathbf{x})]$ be non-zero in $\mathcal{C}$; then it is satisfiable in a model of cardinality $\lambda$ by the completeness theorem for $\mathcal{L}_{\lambda^+, \lambda}$, so that there is $\mathbf{a}$ in $\mathcal{M}$ with $\mathcal{M} \models \psi(\mathbf{a})$. Let $p$ be the type $\{\phi(\mathbf{x}): \mathcal{M} \models \phi(\mathbf{a})\}$.  If $p$ was non-isolated, there would be a model omitting it, i.e., there would exist a model $\mathcal{N}$ of:

$$\mathbb{T} \cup \{``\text{there are $\lambda$ distinct elements"}\} \cup \{\forall \mathbf{x} \bigvee_{\mathcal{M} \models \phi(\mathbf{a})} \neg \phi(\mathbf{x})\}$$
\\
This is impossible since by hypothesis $\mathcal{M}$ and $\mathcal{N}$ are $\mathcal{L}_{\infty, \lambda}$-elementarily equivalent. Therefore, $p$ must be isolated by some $\theta(\mathbf{x})$, which must then be an atom in the Boolean algebra of subobjects of $[\mathbf{x}, \top]$. It follows from this that such algebra is atomic. Moreover, the join of all atoms is the top element (as this is the case in any one, and therefore all, models).

Let us now see that the $\lambda$-classifying topos must be Boolean. Such a topos is built by considering sheaves on $\mathcal{C}$ when equipped with the $\kappa$-topology $\tau$ generated by those jointly epic families of cardinality at most $\lambda$ that corresponds to axioms of (the $\kappa$-geometric Morleyization of) the theory. Let $\mathcal{C}'$ be the full subcategory of $\mathcal{C}$ consisting of non-zero objects, and $\tau'$ the $\kappa$-topology induced by $\tau$. Then the topos $\mathcal{S}h(\mathcal{C}', \tau')$ is still the $\lambda$-classifying topos, but now the $\kappa$-topology $\tau'$ coincides with the coverage $\rho$ consisting of stable nonempty sieves. Indeed, the nontrivial part is showing that a $\rho$-covering sieve $R$ on an object $[\mathbf{x}, \phi]$ of $\mathcal{C}'$ is also $\tau'$-covering. Since $\phi$ is a union of at most $\lambda$ atoms, for each atom there is an arrow in $R$ factoring through it, and since its domain is nonzero, its image must be the whole atom. Choosing one such morphism of $R$ for each atom we get a jointly epic family from $\tau$ contained in $R$. Finally, it follows that the $\lambda$-classifying topos is equivalent to the topos of sheaves on $\mathcal{C}'$ for the double negation topology, which is Boolean.

($\impliedby$) Suppose that the $\lambda$-classifying topos of the theory

$$\mathbb{T} \cup \{``\text{there are $\lambda$ distinct elements"}\}$$
\\
is two-valued and Boolean. Since it is also $\lambda$-separable, it has enough $\lambda$-points (see section \ref{comp}), and in particular it must be atomic. Hence, $\mathcal{C}$ is also atomic. Let $p_i=\{\theta(\mathbf{x_i}): \theta \text{ is an atom in $\mathcal{S}ub([\mathbf{x_i}, \top])$}\}$, where for each $i<\kappa$, $\mathbf{x_i}=x_0x_1...$ up to (but excluding) $i$. Then the family $[\mathbf{x_i}, \theta] \to [\mathbf{x_i}, \top]$ is jointly epic. Therefore, any $\lambda$-point of the topos corresponds to an atomic model of $\mathbb{T}$ of cardinality at least $\lambda$. Since the topos is also two-valued, all such atomic models are $\mathcal{L}_{\lambda^+, \lambda}$-elementarily equivalent. Hence, a back and forth argument shows that any two such models of cardinality $\lambda$ must be $\mathcal{L}_{\infty, \lambda}$-elementarily equivalent.

To complete the proof, notice that atomic toposes are Boolean, while a Boolean topos with enough points must be atomic, and Boolean toposes are two-valued if and only if they are connected.
\end{proof}

Using Theorem \ref{cct} we can now get rid of the Booleanness assumption on $\mathcal{C}$:

\begin{cor}\label{cct2}
Any two models of size $\kappa$ of a $\kappa$-separable topos, whose category of $\kappa$-points has directed colimits, are $\mathcal{L}_{\infty, \kappa}$-elementarily equivalent if and only if it is two-valued and Boolean (alternatively, atomic and connected).
\end{cor}

\begin{proof}
As explained in \cite{espindolad}, the topos $\kappa$-classifies a $\kappa$-coherent theory $\mathbb{T}$. If we let $\mathbb{T}_{B}$ be the theory obtained from $\mathbb{T}$ by adding all instances of excluded middle over $\mathcal{L}_{\kappa^+, \kappa}$, then we have a stable surjection $s: \Sets[\mathbb{T}_B] \twoheadrightarrow \Sets[\mathbb{T}]$ in the category of Grothendieck toposes (indeed, given the pullback $\mathcal{F}$ of $s$ along $\mathcal{G} \to \Sets[\mathbb{T}]$, if $\mathcal{G}$ classifies the geometric theory $\mathbb{S}$, then $\mathcal{F}$ will classify the theory $\mathbb{S}_B$ obtained by adding instances of excluded middle over $\mathcal{L}_{\infty, \omega}$, which is a conservative extension). Since two-valued Boolean toposes are atoms in the lattice of subtoposes of a given topos, it is enough to show that $\Sets[\mathbb{T}]$ has no proper non-degenerate subtoposes if and only if any two models of size $\kappa$ are $\mathcal{L}_{\infty, \kappa}$-elementarily equivalent. Suppose this latter condition holds; consider a subtopos $i: \mathcal{E} \hookrightarrow \Sets[\mathbb{T}]$, and pull it back along $s$. We get a geometric morphism $t: \mathcal{T} \to \Sets[\mathbb{T}_B]$ and a surjection $s': \mathcal{T} \twoheadrightarrow \mathcal{E}$. Now $\Sets[\mathbb{T}]$ has the same models of size $\kappa$ as $\Sets[\mathbb{T}_B]$ has; if we prove that the sequent $\bigwedge_{i<\kappa}\neg \neg \phi_i \vdash_{\mathbf{x}} \neg \neg \bigwedge_{i<\kappa}\phi_i$ holds in $(\mathcal{C}_{\mathbb{T}_B})_{\kappa^+}$, we could run the same proof of Theorem \ref{cct} (since we can apply the omitting types theorem) to show that $\Sets[\mathbb{T}_B]$ is Boolean and two-valued, equivalent in turn to having no proper non-degenerate subtoposes. Therefore, by considering the surjection-embedding factorization of $t$, this is equivalent to either $t$ being a surjection or $\mathcal{T}$ being degenerate. In the first case, it follows that the composite $st$ is a surjection, and since $st \simeq is'$, that $i$ must be a surjection, in which case $\mathcal{E}$ is equivalent to $\Sets[\mathbb{T}]$. In the second case, since $s'$ is a surjection, it follows that $\mathcal{E}$ must be degenerate. 

  It remains to prove that $\bigwedge_{i<\kappa}\neg \neg \phi_i \vdash_{\mathbf{x}} \neg \neg \bigwedge_{i<\kappa}\phi_i$ holds in $(\mathcal{C}_{\mathbb{T}_B})_{\kappa^+}$, for which it is enough to show that the stable surjection $s: \Sets[\mathbb{T}_B]_{\kappa^+} \twoheadrightarrow \Sets[\mathbb{T}]_{\kappa^+}$ is actually an equivalence. Now, the pullback of the double negation subtopos $\mathcal{S}$ of $\Sets[\mathbb{T}]_{\kappa^+}$ along $s$ is equivalent to $\mathcal{S}$.\footnote{If the site for $\mathcal{S}$ is $\mathcal{C}_{\mathcal{S}}$, then the site of the pullback is the free $\kappa$-Boolean extension of $\mathcal{C}_{\mathcal{S}}$, which coincides with $\mathcal{C}_{\mathcal{S}}$ if this latter is already $\kappa$-Boolean.} Since the embedding into $\Sets[\mathbb{T}_B]_{\kappa^+}$ must be dense as this latter is two-valued, we get that $\Sets[\mathbb{T}_B]_{\kappa^+}$ and $\Sets[\mathbb{T}]_{\kappa^+}$ have equivalent double negation subtoposes. This makes $s$ an embedding and thus an equivalence, which completes the proof. 
\end{proof}

 Recall from \cite{rosicky} that an object $M$ of an accessible category is said to be $\lambda$-closed if all morphisms $M \to N$ are $\lambda$-pure. In model-theoretic terms, this means that all morphisms $M \to N$ reflect the truth of $\lambda$-coherent existential formulas with (less than $\lambda$) parameters from $M$ whenever the category is axiomatizable in $\lambda$-coherent logic. One consequence of Corollary \ref{cct2} that will be of use is the following:

\begin{lemma}\label{inj}
Assume $\lambda^{<\lambda}=\lambda$. If $\mathcal{K}$ is categorical in $\lambda > \kappa$, the subcategory of $\mathcal{K}_{\geq \lambda}$ where we just consider the monomorphisms is closed under directed colimits.
\end{lemma}

\begin{proof}
Let $\mathbb{T}$ be the $\mathcal{L}_{\kappa, \theta}$ theory axiomatizing $\mathcal{K}$. The $\lambda$-classifying topos of $\mathbb{T}_{\lambda}$ is Boolean by Corollary \ref{cct2}, so that the evaluation functor $\mathcal{C}_{\mathbb{T}_{\lambda}} \to \Sets^{\mathcal{K}_{\geq \lambda}}$ is Boolean and the sequent $\top \vdash_{\mathbf{x}} \phi \vee \neg \phi$ is forced in every model above the categoricity cardinal. In particular, all morphisms in the restriction of $\mathcal{K}_{\lambda}$ to the monomorphisms will be furthermore $\omega$-pure, and since directed colimits of $\omega$-pure morphisms are $\omega$-pure, such restriction is closed under directed colimits.
\end{proof}

  Continuing with the same notation, we have now:

\begin{lemma}\label{at}
Assume $GCH$. If $\lambda$ is a categoricity cardinal, then the model of size $\lambda$ is $\lambda$-closed.
\end{lemma}

\begin{proof} 
Assume first that $\lambda$ is regular. Then the model $M$ of size $\lambda$ is $\lambda$-atomic. Given a tuple $\mathbf{c}$ and an embedding $M \to N$ with $N \models \exists \mathbf{x} \phi(\mathbf{x}, \mathbf{c})$, let $\phi_0$ be the complete formula that isolates the type of $\mathbf{c}$ in $M$. Then, since $\phi_0(\mathbf{y}) \wedge \exists \mathbf{x} \phi(\mathbf{x}, \mathbf{y})$ is consistent with $\mathbb{T}_{\lambda}$, we must have $\phi_0(\mathbf{y}) \vdash_{\mathbf{y}} \exists \mathbf{x} \phi(\mathbf{x}, \mathbf{y})$ within $\mathbb{T}_{\lambda}$, and in particular $M \models \exists \mathbf{x} \phi(\mathbf{x}, \mathbf{c})$.

In the case when $\lambda$ is limit, we consider first the topos $\Sets[\mathbb{T}_{\lambda}^B]_{\lambda}$, where $\mathbb{T}_{\lambda}^B$ is obtained from $\mathbb{T}_{\lambda}$ by adding all instances of excluded middle for $<\lambda$-coherent formulas. Apply now the same proof as in Theorem \ref{ott} to the $\lambda$-type $\{\phi_i(\mathbf{x})\}_{i<\lambda}$ realized by $\mathbf{c}$ in $M$ to show that, if it is not isolated, there is a Grothendieck coverage $\tau$ on $(\mathcal{C}_{\mathbb{T}_{\lambda}^B})_{\lambda} := \lim \{(\mathcal{C}_{\mathbb{T}_{\kappa}^B})_{\kappa}: \kappa<\lambda \}$ which contains covers corresponding to the omission of the type. Now consider the natural functor $((\mathcal{C}_{\mathbb{T}_{\lambda}^B})_{\lambda}, \tau) \to ((\mathcal{C}_{\mathbb{T}_{\lambda}^B})_{\lambda^+}, \tau)$ inducing a morphism of the corresponding toposes. The codomain topos has enough $\lambda^+$-points, and  these have directed colimits (as can be seen by noticing that the stable surjection $s: \Sets[\mathbb{T}_{\lambda}^B]_{\lambda^+} \twoheadrightarrow \Sets[\mathbb{T}_{\lambda}]_{\lambda^+}$ is actually an equivalence, proved analogously to the last paragraph of the proof of Corollary \ref{cct2}), so the type in question must be omitted by some model of size at most $\lambda^+$. But this model must contain the model of size $\lambda$ as a submodel which realizes the type, which is absurd. This proves that the type is isolated by a complete Boolean formula and hence that $\Sets[\mathbb{T}_{\lambda}^B]_{\lambda}$ is two-valued and Boolean. With the same argument as in Corollary \ref{cct2}, we can conclude that also $\Sets[\mathbb{T}_{\lambda}]_{\lambda}$ is two-valued and Boolean, whence atomic (as the point of size $\lambda$ must be a surjection) and then we can proceed as in the regular case.
\end{proof}

\begin{rmk}\label{catib}
It follows with the same proof idea of Lemma \ref{at} that if $\mathbb{T}$ is $\lambda$-categorical then $\Sets[\mathbb{T}_{\lambda}]_{\lambda}$ is two-valued and Boolean, even if $\lambda$ is singular (in which case by $\Sets[\mathbb{T}_{\lambda}]_{\lambda}$ we understand the limit of the toposes $\Sets[\mathbb{T}_{\nu}]_{\nu}$ for regular $\nu<\lambda$).
\end{rmk}

We get now:

\begin{lemma}\label{lc}
Assume $GCH$ and amalgamation. Then categoricity cardinals are closed in the class of all cardinals.
\end{lemma}

\begin{proof} 
Let $\lambda$ be a limit of categoricity cardinals. Any model of size $\lambda$ is a directed colimit of models of size $\nu$ for $\nu$ a categoricity cardinal. By Lemma \ref{at}, each of the models in the colimit is $\nu$-closed. Since the L\"owenheim-Skolem theorem (namely, that every subset $A \subseteq M$ of a model of size $\kappa^+$ with $|A|=\kappa$ contains a submodel $N \supseteq A$ of size $\kappa$) still holds in terms of internal size, and under $SCH$, for those cardinals $\kappa$ which have cofinality less than $\theta$, it follows that the model of size $\lambda$ is $\lambda$-closed. Since amalgamation holds, it is also $\lambda$-saturated, and there is only one such model of size $\lambda$ up to isomorphism (cf. \cite{rosicky}).
\end{proof}

\section{Classifying toposes for saturated models}

Recall from \cite{rosicky} the categorical definition of $\kappa$-saturated model as a model $M$ such that for every $p: N \to N'$ with $N, N'$ $\kappa$-presentable, each $q: N \to M$ extends to some $q': N' \to M$:
\begin{center}
\begin{tikzcd}
N' \arrow[rrd, "q'", dashed] &  &   \\
N \arrow[rr, "q"',] \arrow[u, "p",] &  & M
\end{tikzcd}
\end{center}
\noindent
This notion coincides with the usual model-theoretic notion in finitary theories, and can be used in the infinitary case to simplify considerably the arguments of \cite{shma}.

  We are now going to prove the following:

\begin{thm}\label{sat}
Assume $GCH$ and amalgamation. Then the $\kappa^+$-classifying topos $\Sets[\mathbb{T}^{sat}_{\kappa^+}]_{\kappa^+}$ of the theory of $\kappa^+$-saturated models is precisely $\mathcal{S}h(\mathcal{K}_{\kappa}^{op}, \tau_D)$, where $\tau_D$ is the dense topology (i.e., it is the double negation subtopos of $\Sets^{\mathcal{K}_{\kappa}}$).
\end{thm}

\begin{proof}
Consider the following diagram:

\begin{center}
\begin{tikzcd}
(\mathcal{C}_{\mathbb{T}_{\kappa}})_{\kappa^+} \arrow[rr, "ev"] &  & \Sets^{\mathcal{K}_{\kappa}} \arrow[dd, "f^*"'] \arrow[rrd, "M \cong \varinjlim_iev_{N_i}"] &  &       \\
                                                                &  &                                                                                             &  & \Sets \\
                                                                &  & {\mathcal{S}h(\mathcal{K}_{\kappa}^{op}, \tau_D)} \arrow[rru, dashed]                       &  &      
\end{tikzcd}
\end{center}
\noindent
The existence of directed colimits readily implies that the dense topology in $(\mathcal{C}_{\mathbb{T}_{\kappa}})_{\kappa^+}$ is a $\kappa^+$-topology (transfinite composites up to $\kappa$ of dense covers are dense). Whence, the double negation subtopos is a $\kappa^+$-topos and the sheafification functor preserves $\kappa^+$-small colimits. Since there is amalgamation, we have that a model $M$ is $\kappa^+$-saturated if and only if for every $p: N \to N'$ in $\mathcal{K}_{\kappa}$, each $q: N \to M$ extends to some $q': N' \to M$ (see the diagram above). This is the same as saying that $M: \Sets^{\mathcal{K}_{\kappa}} \to \Sets$ maps $p^*: [N', -] \to [N, -]$ to an epimorphism, since, writing $M \cong \varinjlim_iev_{N_i}$ as a $\kappa^+$-filtered colimit of evaluations, and given that any $N$ in $\mathcal{K}_{\kappa}$ is $\kappa^+$-presentable, we have:

$$\varinjlim_iev_{N_i}([N, -]) = \varinjlim_i [N, N_i] \cong [N, \varinjlim_i N_i] \cong [N, M]$$
\noindent
It follows that $M$ is $\kappa^+$-saturated if and only if it factors through $f^*: \Sets^{\mathcal{K}_{\kappa}} \to \mathcal{S}h(\mathcal{K}_{\kappa}^{op}, \tau_D)$, since the amalgamation property in $\mathcal{K}_{\kappa}$ is precisely the right Ore condition on $\mathcal{K}_{\kappa}^{op}$, and so it implies that the dense topology coincides with the atomic topology (every nonempty sieve covers). This finishes the proof.
\end{proof}

\begin{rmk}\label{slc}
The description of the $< \lambda$-classifying topos for $<\lambda$-saturated models can be given, with a similar proof as that of Theorem \ref{sat}, and with the same hypothesis, as the double negation subtopos of the topos of presheaves on models of size less than $\lambda$, which will therefore be equivalent to the pseudolimit of the $\kappa$-classifying toposes for $\kappa$-saturated models.
\end{rmk}

\section{Categoricity and amalgamation}

Before proving Grossberg conjecture, we will now show that some of the properties that hold in $AEC's$ also hold in the context of accessible categories with directed colimits. We start with the following:

\begin{lemma}\label{adm}
Assume that $(\kappa^+)^{\kappa} = \kappa^+$. Then $\mathcal{K}_{\kappa}$ has the amalgamation property if and only if $\Sets^{\mathcal{K}_{\kappa}}$ forces the sequent $\top \vdash_{\mathbf{x}} \neg \exists \mathbf{y} \phi(\mathbf{x}, \mathbf{y}) \vee \neg \neg \exists \mathbf{y} \phi(\mathbf{x}, \mathbf{y})$ for every $\kappa^+$-coherent formula $\phi(\mathbf{x}, \mathbf{y})$.
\end{lemma}

\begin{proof}
If $\mathcal{K}_{\kappa}$ has the amalgamation property then $\Sets[\mathbb{T}_{\kappa}]_{\kappa^+} \cong \Sets^{\mathcal{K}_{\kappa}}$ is a De Morgan topos (see \cite{jlnm}), so that in particular it forces the relevant sequent. Conversely, suppose that such a sequent is forced and consider two morphisms $M \to N$ and $M \to S$. If they cannot be amalgamated, using the L\"owenheim-Skolem property, it follows by the completeness theorem (section \ref{comp}) that the $\mathcal{L}_{\kappa^{++}, \kappa^+}$-theory axiomatized by $\mathbb{T}_{\kappa} \cup Diag^+(N) \cup Diag^+(S)$ is inconsistent (where we use the same constants for elements coming from $M$ in the definition of the positive diagrams). Therefore, there is a $\kappa^+$-coherent formula $\phi$ with constant symbols $\mathbf{c}$ only from $M$ such that $\mathbb{T}_{\kappa} \cup Diag^+(N)$ proves $\neg \exists \mathbf{y} \phi(\mathbf{c}, \mathbf{y})$ (if $\kappa \notin S$ we slightly modify the diagrams). Hence, in $\Sets^{\mathcal{K}_{\kappa}}$, $N \Vdash \neg \exists \mathbf{y} \phi(\mathbf{c}, \mathbf{y})$. Since we have $M \Vdash \neg \exists \mathbf{y} \phi(\mathbf{c}, \mathbf{y}) \vee \neg \neg \exists \mathbf{y} \phi(\mathbf{c}, \mathbf{y})$, we must have $M \Vdash \neg \exists \mathbf{y} \phi(\mathbf{c}, \mathbf{y})$. This however contradicts the fact that $S \Vdash \exists \mathbf{y} \phi(\mathbf{c}, \mathbf{y})$.
\end{proof}

  We now prove that the model-theoretic hypothesis of amalgamation for accessible categories can be deduced from categoricity in a high enough cardinal, assuming that there is a proper class of strongly compact cardinals. We have:

\begin{thm}\label{afsc}
Assume $GCH$. Let $\kappa$ be a strongly compact cardinal and let $\mathcal{K}$ be equivalent to the category of models of some $\mathcal{L}_{\kappa, \kappa}$ theory $\mathbb{T}$. If $\mathcal{K}$ is categorical at $\lambda \geq \kappa$, then $\mathcal{K}_{\geq \kappa}$ has the amalgamation property.
\end{thm}

\begin{proof} Since $\lambda$ is a categoricity cardinal, it follows by Lemma \ref{at} that $N_{\lambda}$, the model of size $\lambda$, is $\lambda$-closed, so in particular it is $\kappa$-closed. Thus, we have that for any $\kappa$-coherent $\sigma(\mathbf{x})$, in the topos $\Sets^{\mathcal{K}_{\geq \kappa, \leq \lambda}}$ the model $N_{\lambda}$ forces the sentence $\forall \mathbf{x}(\sigma(\mathbf{x}) \vee \neg \sigma(\mathbf{x}))$, so that $\Sets^{\mathcal{K}_{\geq \kappa, \leq \lambda}}$ forces its double negation. By the intuitionistic completeness theorem of \cite{espindola}, there is a conservative Heyting embedding $ev: (\mathcal{C}_{\mathbb{T}_{\kappa}})_{\kappa} \to \Sets^{\mathcal{K}_{\geq \kappa, \leq \lambda}}$. Since $(\mathcal{C}_{\mathbb{T}_{\kappa}})_{\kappa}$ is two-valued (by categoricity in $\lambda$), $(\mathcal{C}_{\mathbb{T}_{\kappa}})_{\kappa}$ satisfies $\forall \mathbf{x}(\sigma(\mathbf{x}) \vee \neg \sigma(\mathbf{x}))$. This implies that all morphisms in $\mathcal{K}_{\geq \kappa}$ are $\kappa$-pure; amalgamation in $\mathcal{K}_{\geq \kappa}$ follows now as in the finitary case.
\end{proof}

  Note that $SCH$ holds above a strongly compact cardinal, so that it is also enough to assume the latter for our results. In the rest of this section we will work to remove the large cardinal hypothesis.

  Studying eventual categoricity in an accessible category with directed colimits which is categorical in some cardinal is the same as studying it in the $\mu$-AEC obtained by restricting the morphisms to monomorphisms, since by Lemma \ref{inj} this restriction also has directed colimits. The advantage is that we can now consider the functorial expansion (see \cite{bgrlv} for the definition) of the $\mu$-AEC in question defined by Vasey, who called them substructure functorial expansion, in which we add a $\kappa^+$-small arity predicate $P$ whose interpretation in a model $M$ consists of the image of the underlying structure of a model $N$ of size $\kappa$ embedded in $M$ through a morphism in the $\mu$-AEC. This particular expansion, which gives rise to an isomorphic $\mu$-AEC, has the property that morphisms coincide with substructure embeddings. Moreover, its models of size at least $\kappa$ can be axiomatized as follows, extending further the language with the symbol $\subseteq$:

$$\top \vdash_{\mathbf{x}} \exists \mathbf{y}\left(\bigvee_{M_0 \in S}\psi_{M_0}(\mathbf{y}) \wedge \mathbf{x} \subseteq \mathbf{y} \wedge P(\mathbf{y})\right)$$

$$\top \vdash_{\mathbf{x}\mathbf{y}} (\mathbf{x} \subseteq \mathbf{y} \wedge P(\mathbf{x}) \wedge P(\mathbf{y})) \to \bigvee_{(M_0, M_1) \in T}\psi_{(M_0, M_1)}(\mathbf{x}, \mathbf{y})$$

$$\top \vdash_{\mathbf{x}\mathbf{y}} \mathbf{x} \subseteq \mathbf{y} \leftrightarrow \bigwedge_{i \in I}\bigvee_{j \in J}x_i=y_j$$
\noindent
Here $S$ is a skeleton of the subcategory of models of size $\kappa$, $T$ is the set of pairs $(M_0, M_1)$
with a morphism in the $\mu$-AEC and $M_0, M_1 \in S$, while $\psi_{M_0}, \psi_{M_0, M_1}$ are conjunctions of
atomic and negated atomic formulas of the extended language such that $\psi_{M_0}(\mathbf{z})$ holds if and only if $\mathbf{z}$ is isomorphic to $M_0$, and $\psi_{M_0, M_1}(\mathbf{z}, \mathbf{w})$ holds if and only if $(\mathbf{z}, \mathbf{w})$ is isomorphic to $(M_0, M_1)$. With these definitions, it should be clear that the previous sequents really axiomatize the $\mu$-AEC in question: the first sequent is just an axiomatic reformulation of the downward L\"owenheim-Skolem theorem, while the second just expresses that the morphisms are determined solely by the morphisms between models of size $\kappa$. The third sequent is separately added for convenience in what follows.

  If we assume that we have categoricity at $\kappa$, we can now get an axiomatization of an isomorphic $\mu$-AEC which can be entirely rewritten through sequents in the $\kappa^+$-$Reg_{\neg}$ fragment. This is an intuitionistic fragment of first-order logic which contains no disjunctions, obtained from the $\kappa^+$-regular fragment by adding $\bot$, together with the axioms $\bot \vdash_{\mathbf{x}} \phi$ and the axioms for $\neg$ that make it into a negation operator applicable only to atomic formulas. Indeed, in the first sequent above the disjunction reduces to a single disjunct since we have categoricity at $\kappa$, while the second and third sequent above have the general form of universal sentences $\forall \mathbf{z} \bigvee_{i \in I} \bigwedge_{j \in J} \psi_{ij}$, and each such sentence is equivalent to the set of sequents $\{\exists \mathbf{z} \bigwedge_{i \in I} \neg \psi_{if(i)} \vdash \bot\}_{f \in J^I}$.

  The $\kappa^+$-$Reg_{\neg}$ fragment contains the $\kappa^+$-$Reg_{\bot}$ subfragment, not containing the symbol $\neg$. The syntactic category $\mathcal{C}$ of any $\kappa^+$-$Reg_{\neg}$ theory can be studied through the category $\mathcal{K}_{\kappa^+}^r$ of its $\kappa^+$-$Reg_{\bot}$ models (models of the $\kappa^+$-$Reg_{\bot}$ internal theory of $\mathcal{C}$, also known as the $\kappa^+$-$Reg_{\bot}$ Morleyization of the $\kappa^+$-$Reg_{\neg}$ theory). These latter are in particular $\kappa^+$-regular models for the extended signature in which there is an extra propositional symbol $\bot$ and one predicate symbol $S$ for each negated atomic formula $\neg R$ and where the axioms of the theory contain all axioms obtained from formally replacing $\neg R$ by $S$ in each $\kappa^+$-$Reg_{\neg}$ axiom and, additionally, all those axioms of the form $\bot \vdash_{\mathbf{x}} \phi$ and $R \wedge S \vdash_{\mathbf{x}} \bot$.

   If $(\mathcal{C})_{\lambda^+}^r$ is the syntactic category of the $\lambda^+$-$Reg_{\bot}$ theory with the same axioms as the $\kappa^+$-$Reg_{\bot}$ theory of $\mathcal{C}$, then its $\lambda^+$-classifying topos $\mathcal{S}h((\mathcal{C})_{\lambda^+}^r, \tau)$ (where $\tau$ is the $\lambda$-$Reg_{\bot}$ coverage) will be precisely equivalent to the presheaf topos $\Sets^{\mathcal{K}_{\geq \kappa^+, \leq \lambda}^r}$, as can be seen as a special case of Theorem \ref{lk}. In particular, the embedding $(\mathcal{C})_{\lambda^+}^r \to \Sets^{\mathcal{K}_{\geq \kappa^+, \leq \lambda}^r}$ will preserve $\neg$ since it can be identified with Yoneda embedding, which preserves any right adjoint to pullback functors that might exist, see \cite{bj}).

  Using the compactness of $\kappa^+$-$Reg_{\bot}$ logic, it is also easy to verify that the canonical functor $F: \mathcal{C} \to (\mathcal{C})_{\lambda^+}^r$ also preserves $\neg$. For if given a $\lambda^+$-regular sentence $\exists \mathbf{x} \bigwedge_{i<\lambda}\phi_i$ we have $\exists \mathbf{y} \bigwedge_{i<\lambda}\phi_i \wedge R \vdash_{\mathbf{x}} \bot$ in $\lambda^+$-$Reg_{\bot}$ logic, there must be a $\kappa^+$-regular sentence $\exists_{i \in T} y_i \bigwedge_{i \in T}\phi_i$, for some subset $T \subset \lambda$ of size at most $\kappa$, such that $\exists_{i \in T} y_i \bigwedge_{i \in T}\phi_i \wedge R \vdash_{\mathbf{x}} \bot$ in $\kappa^+$-$Reg_{\bot}$ logic, from which our result follows.

  It follows, in fact, that the evaluation functor $ev: \mathcal{C} \to \Sets^{\mathcal{K}_{\geq \kappa^+, \leq \lambda}^r}$, the composite of Yoneda embedding with $F$, preserves $\neg$,\footnote{It is also possible to give a direct proof of this fact, using the compactness of $\kappa^+$-$Reg_{\bot}$ logic, with the same arguments as in the proof of Joyal's theorem, according to which $ev: \mathcal{C} \to \Sets^{Mod(\mathcal{C})}$ preserves universal quantification when $Mod(\mathcal{C})$ is the category of coherent models of the Heyting category $\mathcal{C}$. This is worked out in the author PhD thesis for the more general disjunction-free fragment.} which in particular means that the interpretation of $S$ in the presheaf topos will be precisely that of $\neg R$.

  Note that, if we add to the $\kappa^+$-$Reg_{\neg}$ axiomatization above all instances of excluded middle for atomic formulas, we get an axiomatization of (an isomorphic copy of) the $\mu$-AEC. This fact will be essential for the proof of the following:

\begin{thm}\label{gc}
(Grossberg conjecture for $\nu$-AEC's) Assume $SCH$ (only if $\nu>\omega$). Categoricity in $\kappa$ and $\mu>\kappa$ implies that $\mathcal{K}_{\geq \mu}$ has amalgamation.
\end{thm}

\begin{proof}
Choose a regular $\lambda \geq \mu$, and suppose first that $(\lambda^+)^{\lambda}=\lambda^+$. Let $\mathbb{T}^r$ be the $\kappa^+$-$Reg_{\bot}$ Morleyization of the $\kappa^+$-$Reg_{\neg}$ axiomatization of $\mathbb{T}_{\kappa^+}$ and $\mathcal{K}_{\geq \kappa^+, \leq \lambda}^r$ its category of $\kappa^+$-$Reg_{\bot}$ models. Then, by the observation above, the subtopos $\Sets^{\mathcal{K}_{\geq \kappa^+, \leq \lambda}} \hookrightarrow \Sets^{\mathcal{K}_{\geq \kappa^+, \leq \lambda}^r}$ is dense, as it $\lambda^+$-classifies the quotient theory obtained by adding axioms of the form $\top \vdash_{\mathbf{x}} R \vee S$, which are instances of excluded middle for atomic formulas $R$, since $[[S]] = [[\neg R]]$ in $\Sets^{\mathcal{K}_{\geq \kappa^+, \leq \lambda}^r}$. Also, the sheafification functor is given by the transpose $i^*$ of the inclusion functor $i: \mathcal{K}_{\geq \kappa^+, \leq \lambda} \to \mathcal{K}_{\geq \kappa^+, \leq \lambda}^r$, as can be verified syntactically, and this will be a $\lambda^+$-geometric morphism since $i^*$ preserves all limits and colimits. Note now that the sheafification functor $j: \Sets^{\mathcal{K}_{\geq \kappa^+, \leq \lambda}^r} \to \mathcal{S}h(\mathcal{K}_{\geq \kappa^+, \leq \lambda}^{op}, \tau_D)$ must send $[N', -] \to [N, -]$ to an epimorphism whenever $N \to N'$ in $\mathcal{K}_{\geq \kappa^+, \leq \lambda}^r$ is $\kappa^+$-pure. Indeed, in that case it can be amalgamated with any morphism $N \to N''$ in $\mathcal{K}_{\geq \kappa^+, \leq \lambda}^r$ using $\kappa^+$-compactness of $\kappa^+$-$Reg_{\bot}$ logic, making $N \to N'$ a cover for the dense topology. Given then a morphism $N \to N'$ in $\mathcal{K}_{\geq \mu, \leq \lambda}$, it must be $\mu$-$Reg_{\neg}$-elementary by categoricity in $\mu$, so that it is in particular a $\kappa^+$-pure embedding in $\mathcal{K}_{\geq \kappa^+, \leq \lambda}^r$. It follows that the sheafification functor $j': \Sets^{\mathcal{K}_{\geq \kappa^+, \leq \lambda}} \to \mathcal{S}h(\mathcal{K}_{\geq \kappa^+, \leq \lambda}^{op}, \tau_D)$ must send $[N', -] \to [N, -]$ to an epimorphism, so that $N \to N'$ generates a cover for the dense topology in $\mathcal{K}_{\geq \kappa^+, \leq \lambda}$, i.e., amalgamation at $|N|$ holds. We postpone for the last section the arguments that will allow us to eliminate the hypothesis $(\lambda^+)^{\lambda}=\lambda^+$ and other uses of $GCH$ in the case of AEC's and replace it with $SCH$ for general $\nu$-AEC's.
\end{proof}

  In the following section we will make use of Theorem \ref{gc} and assume that eventual  amalgamation holds whenever we have categoricity in a high enough cardinal.

\section{Eventual categoricity}  
  We can finally prove the following:
  
\begin{thm}\label{btwn}
Assume $GCH$ and amalgamation. Suppose the accessible category $\mathcal{K}$ has directed colimits and is categorical in some cardinal $\kappa>\chi$ (where $\chi$ is the Hanf number for categoricity). Then it is also categorical at any $\delta > \kappa$. 
\end{thm}

\begin{proof}
Assume first $\delta=\kappa^+$ for a cardinal $\kappa$ of cofinality at least $\theta$; let us see that $\Sets[\mathbb{T}_{\kappa^+}]_{\kappa^+}$ has a unique $\kappa^+$-point of size $\kappa^+$, up to isomorphism. Using Lemma \ref{inj},  we will assume that all morphisms are monomorphisms, since the categoricity spectrum does not change with such restriction. Consider the diagram:

\begin{displaymath}
\xymatrix{
\Sets^{\mathcal{K}_{\kappa}} \ar@{ ->}[rrdd]^{N_{\lambda}} \ar@{->}[dd] & & \\
& & \\
\Sets[\mathbb{T}_{\kappa^+}]_{\kappa^+} \ar@{->}[dd] & & \Sets \\
& & \\
\mathcal{S}h(\mathcal{K}_{\kappa}^{op}, \tau_D) \ar@{-->}[rruu] &  &  
}
\end{displaymath}
\noindent
for a model $N_{\lambda}$ of size $\lambda> \kappa^+$ that is a categoricity cardinal. Note that $\mathcal{S}h(\mathcal{K}_{\kappa}^{op}, \tau_D)$ is a $\kappa^+$-topos since $\tau_D$ is a $\kappa^+$-topology (which follows in turn from the fact that $\mathcal{K}_{\kappa}$ has $\kappa^+$-directed colimits). Here $\Sets[\mathbb{T}_{\kappa^+}]_{\kappa^+}$ is a dense subtopos of $\Sets[\mathbb{T}_{\kappa}]_{\kappa^+} \cong \Sets^{\mathcal{K}_{\kappa}}$ since any model in $\mathcal{K}_{\kappa}$ embeds in a model of size $\kappa^+$ and thus every non-zero $\kappa^+$-coherent sentence in $\Sets^{\mathcal{K}_{\kappa}}$ must be non-zero in $\Sets[\mathbb{T}_{\kappa^+}]_{\kappa^+}$. It follows that the double negation subtopos factors through $\Sets[\mathbb{T}_{\kappa^+}]_{\kappa^+}$. Moreover, it is in fact a subtopos of $\Sets[\mathbb{T}_{\kappa^+}]_{\kappa^+}$ and the sheafification functor (which is the same as the embedding into $\Sets^{\mathcal{K}_{\kappa}}$ followed by the sheafification of this latter) must preserve $\kappa^+$-small limits.

  Now $N_{\lambda}$ is $\lambda$-closed, whence it is $\kappa^+$-closed, and therefore $\kappa^+$-saturated, inducing the bottom morphism above. We will prove that any sequent $\vartheta \vdash_{\mathbf{x}} \psi$ valid in $\mathcal{S}h(\mathcal{K}_{\kappa}^{op}, \tau_D)$ is already valid in $\Sets[\mathbb{T}_{\kappa^+}]_{\kappa^+}$. It would then follow that $\Sets[\mathbb{T}_{\kappa^+}]_{\kappa^+} \cong \mathcal{S}h(\mathcal{K}_{\kappa}^{op}, \tau_D)$ and every model of size $\kappa^+$ is $\kappa^+$-saturated. Categoricity in $\kappa^+$ follows (see \cite{rosicky}). 

  In case $\kappa$ has cofinality less than $\theta$, since models in $\mathcal{K}_{\kappa}$ could have cardinality $\kappa^+$, the topos $\Sets[\mathbb{T}_{\kappa^+}]_{\kappa^+}$ must be replaced by a topos whose $\kappa^+$-points are precisely those models of internal size at least $\kappa^+$. We can take instead the topos of equivariant sheaves on the topological groupoid of models of internal size $\kappa^+$ (cf. section \ref{iaf}), which will be precisely the $\kappa^+$-classifying topos of their theory $\mathbb{T}_{\kappa^+}^{int}$. Indeed, if that was not the case, the topos would have some $\kappa^+$ point of $\Sets^{\mathcal{K}_{\kappa}}$ which is not of internal size $\kappa^+$, and the only possibility is that it has the model $M$ of internal size $\kappa$ as the extra $\kappa^+$-point, in which case it would coincide with $\Sets^{\mathcal{K}_{\kappa}}$, which would then have a conservative embedding into a topos with a jointly conservative set of models of internal size $\kappa^+$. Let us see that this implies that $M$ has a $\kappa^+$-pure embedding into a model of internal size $\kappa^+$ (which is impossible, of course). Let $\phi$ be the conjunction of the diagram of the subset $\mathbf{c}$ of cardinality $\kappa$ in $M$ which is a concrete directed colimit of models of smaller size, plus the negation of each existential sentence with parameters from $\mathbf{c}$ which does not hold in $M$. If no model of size $\kappa^+$ satisfied $\phi$, the sequent $\phi \vdash_{\mathbf{x}} \bot$ would be provable in $\mathbb{T}_{\kappa^+}^{int}$, whence by conservativity it would also be provable in $\Sets^{\mathcal{K}_{\kappa}}$, which is not possible since $M \Vdash \phi(\mathbf{c})$. In summary, the topos of equivariant sheaves on the topological groupoid of models of internal size $\kappa^+$ is precisely the $\kappa^+$-classifying topos of their theory $\mathbb{T}_{\kappa^+}^{int}$. Moreover, $\mathbb{T}_{\kappa^+}^{int}$ is a quotient of $\mathbb{T}_{\kappa^+}$, and thus this sheaf topos is a subtopos of $\Sets^{\mathcal{K}_{\kappa}}$ and it is easy to check that the embedding is dense. Therefore, we can continue with the proof above replacing $\mathbb{T}_{\kappa^+}$ with $\mathbb{T}_{\kappa^+}^{int}$.   

  Now consider the $\kappa^+$-coherent theory whose axioms are those of $\mathbb{T}_{\kappa^+}$. By the completeness theorem of section \ref{comp}, this theory admits a jointly conservative set of $\kappa^+$-coherent prime theories in a new language where we add $\kappa^+$ many new constants, each of which is obtained from $\mathbb{T}_{\kappa^+}$ by adding $\kappa^+$-coherent sentences as axioms. We build now a poset of such prime theories in extended languages intended to form a Kripke model, as explained in section {\ref{iil}}. We make however a modification to the construction: for each prime theory $\Gamma$ over a set $\mathbf{c}$ of $\kappa^+$ many constants we consider the theory $\Gamma' := \Gamma \cup \{d=c_i \vdash \bot\}_{i<\kappa^+}$ where $d$ is a new constant, and for this latter we compute the conservative prime theories over the extended language and define them as the theories at the successor nodes.

 At the level of the syntactic categories, $\mathcal{C}_{\Gamma'}$ is the (pseudo-)colimit of $\mathcal{C}_{\Gamma_j} := \mathcal{C}_{\Gamma}/(\{\}, \bigwedge_{i<j}d \neq c_i)$ in the $2$-category of $\kappa^+$-coherent theories and $\kappa^+$-coherent morphisms. Since this colimit is in fact a $\kappa^+$-Heyting category, it will be equivalent to the (pseudo-)colimit computed in the $2$-category of $\kappa^+$-Heyting theories and $\kappa^+$-Heyting morphisms.

  At limit ordinals, we consider the union of the theories corresponding to each branch; this might not be prime, but if consistent, we extend it to prime theories to form the nodes of the next level of the poset. For the resulting frame to be a Kripke model, we need to prove that $\Gamma'$ is conservative over $\Gamma$. For this, in turn, it is enough to show that any $\kappa^+$-coherent model of $\mathbb{T}_{\kappa^+}$ has a proper $\kappa^+$-pure extension. We claim that, moreover, we can always find a proper $\kappa^+$-Boolean extension of any model $M$ of internal size $\kappa^+$. To prove it, suppose otherwise. Then the topos $\Sets[\mathbb{T}_{\kappa^+}^B]_{\kappa^{++}}/[M, -] \cong \Sets^{\mathcal{K}_{\kappa^+}^B}/[M, -] \cong \Sets^{M/\mathcal{K}_{\kappa^+}^B}$, (where $\mathcal{K}_{\kappa^+}^B$ consists of the models in $\mathcal{K}_{\kappa^+}$ and all its $\kappa^+$-Boolean homomorphisms, and where $\mathbb{T}_{\kappa^+}^B$ is $\mathbb{T}_{\kappa^+}$ plus all those instances of excluded middle  for $\kappa^+$-coherent formulas), would be two-valued and Boolean. On the other hand, we have a stable surjection $\Sets^{M/\mathcal{K}_{\kappa^+}^B} \twoheadrightarrow \Sets^{M/\mathcal{K}_{\kappa^+}}$; this can be seen by considering first the stable surjection $\Sets^{\mathcal{K}_{\kappa^+}^B} \cong  \Sets[\mathbb{T}_{\kappa^+}^B]_{\kappa^{++}} \twoheadrightarrow \Sets[\mathbb{T}_{\kappa^+}]_{\kappa^{++}} \cong \Sets^{\mathcal{K}_{\kappa^+}}$ (which, at the level of the theories that they $\kappa^{++}$-classify, just adds instances of excluded middle for $\kappa^+$-coherent formulas). Then we consider the pullback functor to the slice $\Sets^{\mathcal{K}_{\kappa^+}} \to \Sets[\mathbb{T}_{\kappa^+}]_{\kappa^{++}}/[M, -] \cong \Sets^{M/\mathcal{K}_{\kappa^+}}$, which is a geometric morphism along whose direct image we take the following (pseudo-)pullback:

\begin{center}
\begin{tikzcd}
\Sets^{\mathcal{K}_{\kappa^+}^B} \arrow[rr, two heads]                                 &  & \Sets^{\mathcal{K}_{\kappa^+}}                                 \\
                                                                                       &  &                                                                \\
{\Sets[\mathbb{T}_{\kappa^+}^B]_{\kappa^{++}}/[M, -]} \arrow[uu] \arrow[rr, two heads] &  & {\Sets[\mathbb{T}_{\kappa^+}]_{\kappa^{++}}/[M, -]} \arrow[uu]
\end{tikzcd}
\end{center}

Then the (pseudo-)pullback is precisely $\Sets[\mathbb{T}_{\kappa^+}^B]_{\kappa^{++}}/[M, -] \cong \Sets^{M/\mathcal{K}_{\kappa^+}^B}$, as can be verified using the universal property of the slice. This proves that indeed $\Sets^{M/\mathcal{K}_{\kappa^+}^B} \twoheadrightarrow \Sets^{M/\mathcal{K}_{\kappa^+}}$ is a stable surjection. Since being two-valued and Boolean is equivalent to having no proper non-degenerate subtoposes, by an argument analogous to that of the proof of Corollary \ref{cct2}, we would conclude that, since $\Sets^{M/\mathcal{K}_{\kappa^+}^B}$ is two-valued and Boolean, so would $\Sets^{M/\mathcal{K}_{\kappa^+}}$ be, which is absurd since $M$ is not maximal (which is in turn a consequence of categoricity and amalgamation). This finishes the proof of our claim.

  Next, we claim that in the poset of theories there is a branch of height $\kappa^{++}$ consisting of consistent theories whose union, in the logic $\mathcal{L}_{\kappa^{+++}, \kappa^{++}}$, is still consistent. To see this, note first that the corresponding poset of prime theories can be presented as a Kripke model of the internal $\kappa^+$-Heyting theory\footnote{When defining the fragment, we make sure to include all $\kappa^+$-coherent formulas which arise in the $\kappa^+$-classifying topos as instances of universal quantification, so that the syntactic category is $\kappa^+$-Heyting. Then we can define its internal $\kappa^+$-Heyting theory as the set of all $\kappa^+$-Heyting sequents, in the same signature, which are true in the category.} of the category $(\mathcal{C}_{\mathbb{T}_{\kappa^+}})_{\kappa^+}$ (where the underlying set of each node consists of the constants in the prime theory, cf. section \ref{iil}). This follows from the property that a prime theory proves (from axioms in the internal $\kappa^+$-Heyting theory) any $\kappa^+$-Heyting formula evaluated in constants $\mathbf{c}$ if and only if that formula is forced at $\mathbf{c}$ by the corresponding node of the tree. In turn, this property follows by a straightforward induction on the complexity of that formula using that the successor theories are jointly conservative over a given node in the extended language, as in section \ref{iil}. Now each prime theory is determined by $\mathbb{T}_{\kappa^+}$ and the conjunction of the rest of its axioms $\phi_f(\mathbf{c})$ (we can assume these axioms are atomic sentences using the primeness of the theory). We claim that $\phi_{f} \vdash_{\mathbf{y}_{f}} \neg \neg \bigvee_{g \in \gamma^{\beta+1}, g|_{\beta}=f} \exists \mathbf{x}_{g} \phi_{g}$ in $(\mathcal{C}_{\mathbb{T}_{\kappa^+}})_{\kappa^{++}}$. To see this, note that for any model of $\phi_{f}(\mathbf{c})$ and any embedding into some other model, this latter is not maximal, and thus has an embedding into a model of $\phi_{f}(\mathbf{c}) \wedge \bigwedge_i d \neq c_i$; then, a straightforward L\"owenheim-Skolem argument (cf. section \ref{lst}) shows that such model satisfying $\phi_{f}(\mathbf{c}) \wedge \bigwedge_i d \neq c_i$ for some $\mathbf{c}$ must satisfy $\phi_g(\mathbf{d})$ for some $g$ and $\mathbf{d}$, i.e., it satisfies $\bigvee_{g \in \gamma^{\beta+1}, g|_{\beta}=f} \exists \mathbf{x}_{g} \phi_{g}$. Whence, in $\Sets^{\mathcal{K}_{\kappa^+}}$ we see that any model forcing $\phi_{f}(\mathbf{c})$ must force $\neg \neg \bigvee_{g \in \gamma^{\beta+1}, g|_{\beta}=f} \exists \mathbf{x}_{g} \phi_{g}$. This shows that indeed the premise $\phi_{f} \vdash_{\mathbf{y}_{f}} \neg \neg \bigvee_{g \in \gamma^{\beta+1}, g|_{\beta}=f} \exists \mathbf{x}_{g} \phi_{g}$ (that of the transfinite transitivity rule up to double negation) holds in $(\mathcal{C}_{\mathbb{T}_{\kappa^+}})_{\kappa^{++}}$. Since the premises in the transfinite transitivity rule for limit ordinals, up to double negation, also hold, trivially, the conclusion of the rule tells us that the theories of the minimal nodes of any bar are, up to a double negation, jointly conservative over the theory at the root (we are using that the sequent $\bigwedge_{i<\kappa^+}\neg \neg \phi_i \vdash_{\mathbf{x}} \neg \neg \bigwedge_{i<\kappa^+} \phi_i$ holds in $(\mathcal{C}_{\mathbb{T}_{\kappa^+}})_{\kappa^{++}}$, a consequence of $\mathcal{K}_{\kappa^+}$ having directed colimits). This proves that a branch of height $\kappa^{++}$ has to exist, since otherwise there would be a bar composed of inconsistent theories, while $\phi_{\emptyset}$ is consistent.
  
  By repeating the argument with the consistent $\kappa^{++}$-coherent theory corresponding to a node of height $\kappa^{++}$ (but now in the logic $\mathcal{L}_{\kappa^{+++}, \kappa^{++}}$), we can build a new Kripke model up to $\kappa^{+++}$, and so on. At limit cardinals $\rho=\sup_n \rho_n$ (with $\rho_n$ successor cardinals), note that $\Sets[\mathbb{T}_{\kappa^+}]_{\rho} := \lim \Sets[\mathbb{T}_{\kappa^+}]_{\rho_n} \cong \lim \Sets^{\mathcal{K}_{\geq \kappa^+, \leq \rho_{n-1}}}$ and a similar computation to the one above shows that consistent theories of cardinality $\rho$ do exist. We perform for them a similar construction, since we can adapt the completeness theorem (as mentioned at the end of section \ref{comp}) for those theories in $\mathcal{L}_{\rho^+, \rho}$ when the axioms and the conclusions are in $\rho$-fragments of $\mathcal{L}_{\rho^+, \rho}$ (which are, by definition, unions of all $\rho_n$-fragments); cf. the last paragraph of section \ref{comp}. The corresponding conservative extensions of each model of size $\rho$ follows now from the fact that any model of size $\rho$ has a proper $<\rho$-pure extension, as can be proven similarly to the arguments above.

  Eventually, we reach a consistent theory of cardinality $\lambda$, and we can construct its corresponding Kripke model. Now the theory at the root proves that there are $\lambda$ distinct elements. Since the model of $\mathbb{T}_{\kappa^+}$ of size $\lambda$ satisfies $\forall \mathbf{x} (\vartheta \to \psi)$, this sentence will be forced in the Kripke model, and so will be proved by the corresponding prime theory at the root node. This shows that $\neg \forall \mathbf{x} (\vartheta \to \psi)$ cannot be $1$ in $(\mathcal{C}_{\mathbb{T}_{\kappa^+}})_{\kappa^+}$, as otherwise all successive consistent theories of each cardinality over which we constructed the Kripke models would prove $\neg \forall \mathbf{x} (\vartheta \to \psi)$, in particular the prime theory at the root node of the latest Kripke model. Therefore, as $(\mathcal{C}_{\mathbb{T}_{\kappa^+}})_{\kappa^+}$ is two-valued, $\forall \mathbf{x} (\vartheta \to \psi)$ has to be $1$ there, as we wanted to show.

  As a last point, note that the same proof above which allowed us to conclude categoricity in $\kappa^+$ from categoricity in $\kappa$ also let us prove categoricity in $\kappa^{++}$ and so on, by an inductive procedure. Categoricity in a $\delta$ which is a limit cardinal is easily handled knowing that $\mathcal{K}$ will be categorical at all $\gamma_i$ for a cofinal sequence of successors $\kappa<\gamma_i<\delta$. Indeed, since $\mathcal{K}$ is $\gamma_i$-categorical, the model of size $\gamma_i$ is $\gamma_i$-saturated, which allows us to successively find a set of compatible isomorphisms between submodels of any two models of size $\delta$ (using directed colimits at limit steps), proving that they are indeed isomorphic (see \cite{rosicky}).
\end{proof}

\begin{rmk}\label{rmkesm}
The Kripke models built during the proof of Theorem \ref{btwn} allow to find the parallel with Proposition 4.27 in \cite{shma}, established there with the use of stability theory. More specifically, if we have categoricity in $\kappa$ but not in $\kappa^+$ (i.e., there is a $\kappa$-saturated model which is not $\kappa^+$-saturated), we have that some formula $\neg \forall \mathbf{x} (\vartheta \to \psi)$ is forced at the root of every Kripke model whose root node has cardinality $\mu$; whence we can conclude that there is a non-saturated model in each cardinality $\mu>\kappa^+$.
\end{rmk}

We will prove in the next section how to replace $GCH$ with $SCH$ in the general case of accessible categories with directed colimits, and to remove it completely in the case of AEC's. Our main result, in view of Theorem \ref{btwn}, is the following:

\begin{thm}\label{main}
(Shelah's eventual categoricity conjecture for accessible categories with directed colimits). Assume $SCH$, and let $\mathcal{K}$ be an accessible category with directed colimits. Then there exists a cardinal $\mu_0$ such that if $\mathcal{K}$ is categorical in some $\lambda \geq \mu_0$, it is categorical in all $\lambda' \geq \mu_0$.
\end{thm}

\begin{proof}
It is enough to take $\mu_0$ the maximum of the Hanf numbers for categoricity and non-categoricity.
\end{proof}

\begin{cor}\label{im}
(Morley's categoricity theorem for infinitary theories) Assume $SCH$, and let $\mathbb{T}$ be a $\mathcal{L}_{\kappa, \theta}$-theory whose category of models has directed colimits. Then there exists a cardinal $\mu_0$ such that if $\mathbb{T}$ is categorical in some $\lambda \geq \mu_0$ in $S$, it is categorical in all $\lambda' \geq \mu_0$ in $S$.
\end{cor}

\begin{exm}\label{exm}
In the case in which we consider the cardinality of the underlying model instead of the internal size, the exceptions of Corollary \ref{im} are necessary. Indeed, generalizing a result of \cite{rlv}, the category of $\mu$-Hilbert spaces is defined as follows. Consider a $\mu$-field $\mathbf{R}$, that is a field of hyperreals containing all ordinals up to $\mu$. The construction of such a $\mu$-field proceeds with the following steps:

\begin{itemize}
\item Take the initial segment of the ordinals up to $\mu$. The natural (Hessenberg) sum and product is defined setting $a+b$ (resp. $a.b$) as the maximum order type of a linear order extending the partial order given by the disjoint union (resp. the direct product). They are associative, commutative and the product distributes over the sum. At each following step, the sum and product operations can be defined similarly to the construction of the real numbers.

\item Build the corresponding ring of $\mu$-integers as pairs of ordinals $(a, b)$.

\item Build the field of fractions of that ring.

\item Take the $\mu$-completion of that field considering all $\mu$-Cauchy $\mu$-sequences of fractions.
\end{itemize}

A $\mu$-Hilbert space is then a Hilbert space over the $\mu$-field $\mathbf{R}$. The category of $\mu$-Hilbert spaces and isometries is then axiomatizable in $\mathcal{L}_{\mu^+, \mu^+}$ (e.g., generalizing the axiomatization of Hilbert spaces described in \cite{rlv}). Given an orthonormal base of size $\lambda$, each element of the $\mu$-Hilbert space has at most $\mu$ nonzero coordinates. As a result, the cardinality is of the form $\lambda^{\mu}$. Assuming $GCH$ (or merely $SCH$), we have:

$$\lambda^{\mu}=\left.
\begin{cases}
\lambda & \text{ if $cf(\lambda)>\mu$ and $2^\mu<\lambda$}\\
\lambda^+ &\text{ if $cf(\lambda) \leq \mu$ and $2^\mu<\lambda$}
\end{cases}
\right.$$

It follows from this that there are no models at cardinals of cofinality less than $\mu$, while eventually there is exactly one $\mu$-Hilbert space of cardinality $\lambda$ whenever $\lambda$ is regular but not a successor of a cardinal of cofinality less than $\mu$, and there are two $\mu$-Hilbert spaces (of internal sizes $\lambda$ and $\lambda^+$) if $\lambda$ is such a successor. (On the other hand, it is categorical in every $\lambda$ with respect to internal size).
\end{exm}

In the context of AEC's, we get:

\begin{cor}\label{aec}
(Shelah's eventual categoricity conjecture for AEC's). Let $\mathcal{K}$ be an AEC. Then there exists a cardinal $\mu_0$ such that if $\mathcal{K}$ is categorical in some $\lambda > \mu_0$, it is categorical in all $\lambda' > \mu_0$. 
\end{cor}

\begin{proof}
It is enough to note that internal and external sizes coincide, so we can use the result of Theorem \ref{main}. 
\end{proof}

\section{Removing GCH}\label{ngch}

So far we have proven our results assuming that $GCH$ holds, which in almost every case was needed to guarantee that every regular cardinal $\kappa$ satisfies $\kappa^{<\kappa}=\kappa$. Our method of proof is such that in several cases we can use forcing to make this cardinal equality true by collapsing $\kappa^{<\kappa}$ to $\kappa$ without affecting the main properties of the models of cardinality less than $\kappa$. As a consequence, in many of the results on AEC's, $GCH$ can be removed, while for $\mu$-AEC's it can be downgraded to $SCH$. This general strategy is illustrated in the following:

\begin{thm}\label{forcing}
In all the previous applications, Theorem \ref{lk} holds without any cardinal arithmetic assumptions. 
\end{thm} 

\begin{proof}
The assumptions we have used were $\kappa^{<\kappa}=\kappa$ and $\lambda^{<\lambda}=\lambda$. The first one was there only to guarantee that the cardinality of formulas of the appropriate $\kappa$-fragment of $\mathcal{L}_{\kappa^+, \kappa}$ was less than $\lambda$. Since $\mathcal{K}$ is axiomatizable in $\mathcal{L}_{\kappa, \theta}(\Sigma)$ (for a signature $\Sigma$ with only binary relations), it is enough to require $\kappa^{<\theta}<\lambda$, but as we are using this theorem above the Hanf number for model existence, this is always the case.

 To handle the second assumption, let $\lambda=\gamma^+$. Consider the forcing extension $V[G]$ in which we collapse $\lambda^{<\lambda}$ to $\lambda$. Since this forcing is $<\lambda$-distributive, models of size less than $\lambda$, and their embeddings remain unchanged (we assume they are properly coded by ordinals). Moreover, subsets definable by $\lambda$-geometric formulas of $\mathcal{L}_{\infty, \lambda}$ in $V[G]$ remain definable in $V$. Indeed, any $\lambda$-geometric formula $\phi(\mathbf{x})$ is of the form $\bigvee_{j<2^{\gamma}} \exists_{i<\gamma} x_i \bigwedge_{i<\gamma} \psi_i^j$, where the $\psi_i^j$ are atomic formulas with free variables amongst $\{x_k: k<i\} \cup \mathbf{x}$ and $\mathbf{x}=\cup_{i<\gamma} \mathbf{x_i^j}$ is a union (not necessarily disjoint) of the sets of free variables of each $\psi_i^j$ with $|\mathbf{x_i^j}|=2$. The formula defines a subset $S \subseteq \cup_{i<\gamma}|M|^{\mathbf{x_i^j}}$ in a model $M$ if and only if there is a function $f: S \to V$ whose range is precisely the set of subformulas $\{\exists_{i<\gamma} x_i \bigwedge_{i<\gamma} \psi_i^j: j<2^{\gamma}\}$ such that each $s \in S$ satisfies $f(s)$. Since the model is coded by ordinals (i.e., the underlying set, the functions and relations are all coded by ordinals less than $\gamma$), so can $S$ be coded by ordinals less than $\gamma$, and thus $f$ is in $V$ and the subset is definable in the ground model.

We can now apply the theorem and deduce that the $\lambda$-classifying topos of $\mathbb{T}'$ is equivalent to the presheaf topos $\Sets^{Mod_{\lambda}(\mathbb{T})}$, which reduces to say that the following two conditions are satisfied:

\begin{enumerate}

\item The evaluations $ev_{\phi}$ in $\lambda$-geometric formulas form a generator of $\Sets_{V[G]}^{Mod_{\lambda}(\mathbb{T})}$
\item The evaluation functor $ev: \mathcal{C}_{\mathbb{T}'} \to \Sets_{V[G]}^{Mod_{\lambda}(\mathbb{T})}$ is full on subobjects (here $\mathcal{C}_{\mathbb{T}'}$ is the full syntactic category, not just the one restricted to the $\lambda$-fragment)
\end{enumerate}

  Let us see that both these statements can be expressed as assertions that only involve the existence of functions $f_i: \gamma \to V$; then, by $<\lambda$-distributivity of the forcing, it follows that the functions $f_i$ are actually in $V$, and therefore the statements hold already in the ground model. 

  Condition $1$ boils down to saying that for each functor $F: Mod_{\lambda}(\mathbb{T}) \to \Sets_{V[G]}$ and proper subfunctor $S \rightarrowtail F$ there is a $\lambda$-geometric formula $\phi$ and maps $\eta_M: ev_{\phi}(M) \to F(M)$ satisfying naturality requirements such that $\eta$ does not factor through $S$. In particular, this is true for all functors $F: Mod_{\lambda}(\mathbb{T}) \to \Sets_{V}$. But in that case, since the underlying sets of the models, their functions and relations are all coded by ordinals less than $\gamma$, when $F(M) \subset V$ it follows from $<\lambda$-distributivity that $\eta_M$ is in $V$ as well. Thus, the evaluations in $\lambda$-geometric formulas generate $\mathcal{S}et_{V}^{Mod_{\lambda}(\mathbb{T})}$ and condition $1$ is true in the ground model.

  Condition $2$ is easier to handle, since it is equivalent to saying that any distinguished subset of each $\lambda$-presentable model (which must have size less than $\lambda$) that is preserved by model homomorphisms is definable by a $\lambda$-geometric formula, which is a statement involving definability of fixed subsets of (the underlying sets of) these models, and since those subsets are unchanged by the forcing extension, the statement remains true in the ground model.
\end{proof}

  In analogous ways, several uses of $GCH$ that have been made so far can be eliminated:

\begin{itemize}

\item The use of $GCH$ in Lemma \ref{adm} comes from the equality $(\kappa^+)^{\kappa} = \kappa^+$, whose only purpose was to ensure that the completeness theorem for $\mathcal{L}_{\kappa^{++}, \kappa^+}$ can be applied, at the same time that the $Diag^+(N)$ and $Diag^+(S)$ contain at most $\kappa^{<\theta} \leq \kappa^+$ atomic formulas. To remove the need for $GCH$ in this case, we proceed with the same strategy as in Theorem \ref{forcing}. That is, if the equality does not hold, we consider the forcing extension $V[G]$ in which we collapse $(\kappa^+)^{\kappa}$ to $\kappa^+$. This forcing is $<\kappa^+$-distributive and thus it does not change the category $\mathcal{K}_{\kappa}$. Whence, since we know that $\mathcal{K}_{\kappa}$ satisfies amalgamation in $V[G]$, it already satisfies it in $V$. An analogous argument allows to conclude in Theorem \ref{afsc} amalgamation at $\kappa$ even when $\lambda^{<\lambda}>\lambda$, and to eliminate the uses of $GCH$ in Theorem \ref{gc}. The rest of the uses of $GCH$ in Theorem \ref{btwn}, with one exception that we treat separately below, can be dealt with by using a forcing argument collapsing $(\lambda^+)^{\lambda}$ to $\lambda^+$ to derive the downward categoricity transfer. 

\item Likewise, the use of $GCH$ in Theorem \ref{sat} can be avoided since it was only there to justify the applicability of Theorem \ref{lk}, where it is not needed by Theorem \ref{forcing}. Also, in Lemma \ref{at}, $GCH$ is used for the equality $(\lambda^+)^{\lambda}=\lambda^+$, but note that we can conclude $\lambda$-closedness without $GCH$, by similar forcing arguments (so that Lemma \ref{lc} does not need $GCH$ either). On the other hand, the use of the equality $\kappa^{<\kappa}=\kappa$ in Lemma \ref{inj} was to guarantee the applicability of Corollary \ref{cct2}, where it was crucial to prove the Booleanness of the topos $\Sets[\mathbb{T}_{\kappa}]_{\kappa}$, and to make sure that cardinality coincides with internal size. But this latter condition holds in cardinals in $S$, and the proof of that corollary does not need that equality though to prove that the topos is two-valued, so using this fact we can prove that the morphisms are $\omega$-pure by using the following argument: we know that the topos $\Sets^{\mathcal{K}_{\lambda}}$ is two-valued, and since each model of size $\lambda$ embeds into a $\omega$-closed model of size $\lambda$, $\Sets^{\mathcal{K}_{\lambda}}$ forces $\neg \neg \forall \mathbf{x} (\phi \vee \neg \phi)$ and thus $\forall \mathbf{x} (\phi \vee \neg \phi)$, from which our claim follows.

\end{itemize}

Finally, consider the use of $GCH$ from the proof of Theorem \ref{btwn} for the Kripke completeness theorem. To remove it, define $f(\alpha)=\alpha^{<\alpha}$ and $C(\alpha)=\sup \{f^n(\alpha): n \in \mathbf{N}\}$. Note that $C(\alpha) \times C(\alpha)$ has a well-ordering which is the union of the canonical well-orderings of each $f^n(\alpha) \times f^n(\alpha)$ for $n \geq 1$; this can be used in the completeness proof of section \ref{comp} to find for each $\alpha$-coherent theory a set of prime $C(\alpha)$-coherent theories in a language extended with $C(\alpha)$ many contants and which are jointly conservative over the original theory in the logic $\mathcal{L}_{(C(\alpha))^+, C(\alpha)}$. It suffices then to modify the proof as follows: we define each slice $\mathcal{C}' \cong \mathcal{C}_{\mathbb{T} \cup \phi(\mathbf{c})}$ as before, but now in the definition of the set of jointly epic families of morphisms over each object of the syntactic category, we make sure to put first the set of $\alpha$ many covers $\{C_i \to A\}_{i<\alpha}$ generated, as an $\alpha$-Grothendieck topology, from axioms of the theory, and when considering the set of sections $s_B: B \to A'$ which are used in the transfinite construction of section \ref{comp}, we make sure to put first those arrows, if any, represented by $\alpha$-coherent formulas. We then use a well-ordering $\alpha \times \alpha \to \alpha$ to get as the first $\alpha$ elements of each family $\mathcal{F}(B)$ the pullbacks of the first $\alpha$ coverings along the first $\alpha$ sections. 

  In each of the term models obtained from $\mathcal{C}'$ by the transfinite construction of section \ref{comp} we choose a submodel of size $\alpha$ containing $\mathbf{c}$ (this is possible due to the L\"owenheim-Skolem property). Each of these submodels $M$ determines, by the L\"owenheim-Skolem theorem of section \ref{lst}, a branch $b$ and a further submodel $M_b$ which belongs to the original set of term models, and has cardinality $\alpha$. Let $B_0$ be the set of all such branches, and $T_0$ is the subtree formed by the branches in $B_0$. Since any model of $\mathcal{C}'$ has, as a submodel, some term model of some branch in $B_0$, the morphisms from $\mathcal{C}'$ to the slices over elements of any bar in $T_0$ must be jointly conservative. We claim that, moreover, the set of sections $s_B: B \to A'$ which were used in the transfinite construction have cardinality at most $\alpha$ if $B$ is a node in $T_0$. The reason for this latter claim is that we have restricted ourselves to accessible categories where all morphisms are monomorphisms, so that the formula $\forall xy (x=y \vee x \neq y)$ is forced by any object $B$ of $\mathcal{C}'$; thus, two sections $s_B: B \to A'$ which are different will represent different arrows in the pseudocolimit. 

  Let now start with the syntactic category of the theory in $\alpha$-coherent logic, and take its slice $\mathcal{C}''$. Perform the transfinite construction as before but only up to $\alpha$, using only the $\alpha$ covers $\{C_i \to A\}_{i<\alpha}$ arising from axioms of the theory and the $\alpha$ sections $s_B: B \to A'$ represented by $\alpha$-coherent formulas. It follows that the set of covering families $\mathcal{F}(B)$ over each $B$ in $T_0$ contains, as the first $\alpha$ elements, only those pullbacks, along the $\alpha$ sections $s_B: B \to A'$ represented by $\alpha$-coherent formulas, of the $\alpha$ covers $\{C_i \to A\}_{i<\alpha}$ arising from axioms of the theory. Then, the term models produced by this construction at level $\alpha$ in each branch of $T_0$ will have cardinality at most $\alpha$; likewise, their $\alpha$-coherent theories will also have $\alpha$ many axioms, allowing for the construction to proceed. In this way we can get a jointly conservative set of $\alpha$-coherent prime theories over a language extended with $\alpha$ many constants. This finishes the proof that $GCH$ can be removed in the particular case of AEC's, and downgraded to $SCH$ in general accessible categories with directed colimits.

\begin{rmk}\label{alt}
The argument of the previous paragraphs can be alternatively used to show how to avoid $GCH$ while using the completeness theorem of section \ref{comp} for all the partial results on $\mu$-AEC's, replacing it with $SCH$ and, in the particular case of AEC's, eliminating it completely. As can be seen from the argument, the key is the availability of the downward  L\"owenheim-Skolem property, which allows us to circumvent cardinal arithmetic assumptions. These can be considered, in a sense, to be built in into the L\"owenheim-Skolem property.
\end{rmk}

\bibliographystyle{amsalpha}

\renewcommand{\bibname}{References} % changes the header from Bibliography to References

\bibliography{references}

\newcommand{\etalchar}[1]{$^{#1}$}
\providecommand{\bysame}{\leavevmode\hbox to3em{\hrulefill}\thinspace}
\providecommand{\MR}{\relax\ifhmode\unskip\space\fi MR }
\providecommand{\MRhref}[2]{%
  \href{http://www.ams.org/mathscinet-getitem?mr=#1}{#2}
}
\providecommand{\href}[2]{#2}
\begin{thebibliography}{BGL{\etalchar{+}}16}

\bibitem[AF13]{af}
Steve Awodey and Henrik Forssell, \emph{First-order logical duality}, Annals of
  Pure and Applied Logic \textbf{164} (2013), no.~3, 319--348.

\bibitem[BGL{\etalchar{+}}16]{bgrlv}
Will Boney, Rami Grossberg, Michael Lieberman, Ji\v{r}\'{i} Rosick\'{y}, and
  Sebastien Vasey, \emph{$\mu$-abstract elementary classes and other
  generalizations}, Journal of Pure and Applied Algebra \textbf{220} (2016),
  no.~9, 3048--3066.

\bibitem[BJ98]{bj}
Carsten Butz and Peter Johnstone, \emph{Classifying toposes for first-order
  theories}, Annals of Pure and Applied Logic \textbf{91} (1998), no.~1,
  33--58.

\bibitem[BR12]{br}
Tibor Beke and Jir{\'\i} Rosick{\`y}, \emph{Abstract elementary classes and
  accessible categories}, Annals of Pure and Applied Logic \textbf{163} (2012),
  no.~12, 2008--2017.

\bibitem[Esp19]{espindola}
Christian Esp\'indola, \emph{Infinitary first-order categorical logic}, Annals
  of Pure and Applied Logic \textbf{170} (2019), no.~2, 137--162.

\bibitem[Esp20]{espindolad}
\bysame, \emph{Infinitary generalizations of {D}eligne's completeness theorem},
  The Journal of Symbolic Logic \textbf{85} (2020), no.~3, 1147--1162.

\bibitem[Jec03]{jechst}
Thomas Jech, \emph{Set theory. {T}he {T}hird {M}illenium {E}dition, revised and
  expanded}, Springer Monographs in Mathematics. Springer-Verlag, Berlin
  (2003).

\bibitem[Joh79]{jlnm}
Peter Johnstone, \emph{Conditions related to {D}e {M}organ's law}, Applications
  of sheaves, Springer, 1979, pp.~479--491.

\bibitem[Kar64]{karp}
Carol Karp, \emph{Languages with expressions of infinite length}, North-Holland
  Publishing Company, 1964.

\bibitem[Kel89]{kelly}
Gregory~Maxwell Kelly, \emph{Elementary observations on 2-categorical limits},
  Bulletin of the Australian Mathematical Society \textbf{39} (1989), no.~2,
  301--317.

\bibitem[LRV19]{rlv}
Michael Lieberman, Ji\v{r}\'{i} Rosick\'{y}, and Sebastien Vasey,
  \emph{Internal sizes in $\mu$-abstract elementary classes}, Journal of Pure
  and Applied Algebra \textbf{223} (2019), no.~10, 4560--4582.

\bibitem[Mak90]{makkai}
Michael Makkai, \emph{A theorem on {B}arr-exact categories, with an infinitary
  generalization}, Annals of pure and applied logic \textbf{47} (1990), no.~3,
  225--268.

\bibitem[Moe88]{moerdijk2}
Ieke Moerdijk, \emph{The classifying topos of a continuous groupoid. i},
  Transactions of the American Mathematical Society \textbf{310} (1988), no.~2,
  629--668.

\bibitem[Ros97]{rosicky}
Ji{\v{r}}{\'\i} Rosick{\`y}, \emph{Accessible categories, saturation and
  categoricity}, The Journal of Symbolic Logic \textbf{62} (1997), no.~3,
  891--901.

\bibitem[She09]{shelah}
Saharon Shelah, \emph{Classification theory for abstract elementary classes},
  Studies in Logic: Mathematical logic and foundations, vol. 18, College
  Publications, 2009.

\bibitem[SM90]{shma}
Saharon Shelah and Michael Makkai, \emph{Categoricity of theories in
  $\mathcal{L}_{\kappa, \omega}$, with $\kappa$ a compact cardinal}, Annals of
  Pure and Applied Logic \textbf{47} (1990), no.~1, 41--97.

\bibitem[SS18]{shva}
Shelah Saharon and Vasey Sebastien, \emph{Categoricity and multidimensional
  diagrams}, https://arxiv.org/abs/1805.06291 (2018).

\bibitem[Vas18]{sv3}
Sebastien Vasey, \emph{The categoricity spectrum of large abstract elementary
  classes}, https://arxiv.org/abs/1805.04068 (2018).

\end{thebibliography}

%\begin{thebibliography}{widest entry}

%\end{thebibliography}

\end{document}